\documentclass{amsart}
\usepackage{hyperref}

\newtheorem{theorem}{Theorem}[section]
\newtheorem{lemma}[theorem]{Lemma}
\newtheorem{remark}[theorem]{Remark}
\newtheorem{hypo}[theorem]{Hypothesis {\bf H.}\hspace*{-0.6ex}}

\newcommand{\R}{{\mathbb R}}
\newcommand{\N}{{\mathbb N}}
\newcommand{\Z}{{\mathbb Z}}
\newcommand{\C}{{\mathbb C}}
\newcommand{\M}{{\mathbb M}}
\newcommand{\W}{{\mathbb W}}

\newcommand{\nn}{\nonumber}
\newcommand{\be}{\begin{equation}}
\newcommand{\ee}{\end{equation}}
\newcommand{\bea}{\begin{eqnarray}}
\newcommand{\eea}{\end{eqnarray}}
\newcommand{\ul}{\underline}
\newcommand{\ti}{\tilde}

\newcommand{\spr}[2]{\langle #1 , #2 \rangle}

\newcommand{\I}{\mathrm{i}}
\newcommand{\E}{\mathrm{e}}
\newcommand{\re}{\mathrm{Re}}
\newcommand{\im}{\mathrm{Im}}

\newcommand{\res}{\mathrm{Res}}

\newcommand{\ulz}{\ul{z}}
\newcommand{\hmu}{\hat{\mu}}

\newcommand{\dimuz}{\di_{\ul{\hat{\mu}}}}

\newcommand{\vrc}{\ul{\Xi}_{p_0}}
\newcommand{\hvrc}{\ul{\hat{\Xi}}_{p_0}}

\newcommand{\di}{\mathcal{D}}

\newcommand{\Amap}{\ul{A}_{p_0}}
\newcommand{\amap}{\ul{\alpha}_{p_0}}
\newcommand{\hAmap}{\ul{\hat{A}}_{p_0}}

\newcommand{\Rg}[1]{R_{2g+2}^{1/2}(#1)}

\newcommand{\eps}{\varepsilon}

\newcommand{\sig}{\sigma}
\newcommand{\lam}{\lambda}
\newcommand{\gam}{\gamma}
\newcommand{\om}{\omega}

\numberwithin{equation}{section}

\begin{document}

\title[Scattering Theory for Jacobi Operators]{Scattering Theory for Jacobi Operators with Quasi-Periodic Background}

\author{Iryna Egorova}
\address{Kharkiv National University\\ Ukraine}
\email{\href{mailto:egorova@ilt.kharkov.ua}{egorova@ilt.kharkov.ua}}

\author{Johanna Michor}
\address{Fakult\"at f\"ur Mathematik\\
Nordbergstrasse 15\\ 1090 Wien\\ Austria\\ and International Erwin Schr\"odinger
Institute for Mathematical Physics, Boltzmanngasse 9\\ 1090 Wien\\ Austria}
\email{\href{mailto:Johanna.Michor@esi.ac.at}{Johanna.Michor@esi.ac.at}}

\author{Gerald Teschl}
\address{Fakult\"at f\"ur Mathematik\\
Nordbergstrasse 15\\ 1090 Wien\\ Austria\\ and International Erwin Schr\"odinger
Institute for Mathematical Physics, Boltzmanngasse 9\\ 1090 Wien\\ Austria}
\email{\href{mailto:Gerald.Teschl@univie.ac.at}{Gerald.Teschl@univie.ac.at}}
\urladdr{\href{http://www.mat.univie.ac.at/~gerald/}{http://www.mat.univie.ac.at/\~{}gerald/}}

\thanks{Work supported by the Austrian Science Fund (FWF) under Grant
No.\ P17762, the Austrian Academy of Sciences under DOC-21388, and
INTAS Research Network NeCCA 03-51-6637
}


\begin{abstract}
We develop direct and inverse scattering theory for Jacobi operators which are
short range perturbations of quasi-periodic finite-gap operators. We show
existence of transformation operators, investigate their properties, derive the
corresponding Gel'fand-Levitan-Marchenko equation, and find minimal
scattering data which determine the perturbed operator uniquely.
\end{abstract}

\maketitle

\section{Introduction}

Classical scattering theory deals with the reconstruction of a given Jacobi operator
\begin{equation} \label{eq1}
H u(n)= a(n) u(n+1) + a(n-1) u(n-1) + b(n) u(n) ,
\end{equation}
which is a short range perturbation of the {\em free} one $H_0$ associated with the coefficients
$a(n)=\frac{1}{2}$, $b(n)=0$. This case has been first developed on an informal level by Case in a series
of papers \cite{dinv4} -- \cite{dinv1}. The first rigorous results were established by Guseinov
\cite{gu}, who gave necessary and sufficient conditions for the scattering data to
determine $H$ uniquely under the assumption
\begin{equation} \label{decay}
\sum_n |n| \left( |a(n)- \frac{1}{2}| + |b(n)| \right) < \infty.
\end{equation}
Further extensions were made by Guseinov \cite{gu2}, \cite{gu3}, and
Teschl \cite{tivp}. Additional details and further references can be found, e.g., in \cite{tjac}.

In addition to being of interest on its own, scattering theory can also be used to
solve the initial value problem for the Toda equation via the inverse scattering transform.
This has been formally developed by Flaschka \cite{fl2} (see also \cite{ta} and \cite{fad} for
the case of rapidly decaying sequences) who also worked out the inverse
procedure in the reflection-less case. Further results and an extension of the method to the
entire Toda hierarchy were given by Teschl in \cite{tist} and \cite{tivp}.

The next interesting problem is to replace the free Hamiltonian $H_0$ by one with a periodic
potential. First results in the case of Sturm-Liouville operators have been obtained by Firsova
in a series of papers (see \cite{fir}). For further results, including potentials with different spatial
asymptotics, and additional references see Gesztesy
et al.\ \cite{gnp}. In the discrete case, the investigation has only recently been started by
Boutet de Monvel and Egorova \cite{bdme} and by Volberg and Yuditskii \cite{voyu}, who treat
the case where $H$ has a homogeneous spectrum and is of Szeg\"o class exhaustively
from an operator point of view. Applications to the Toda lattice can be found in Bazargan and
Egorova \cite{baeg} and Boutet de Monvel and Egorova \cite{bdme2}.

\noindent
Finally, let us give a brief overview of the paper:

Section~\ref{secQP} collects some well-known facts from Riemann surfaces and
introduces the necessary notation. Section~\ref{secBA} introduces the Baker-Akhiezer
function and investigates the quasi-momentum map. In the periodic case, where the
integrals can be explicitly computed, this was first done in \cite{perco}. In addition,
we characterize the second solution at the band edges. In Section~\ref{secJS} we
prove existence of Jost solutions and use them to characterize the spectrum of
the perturbed operator. In the periodic case, existence of Jost solutions was first shown
by Geronimo and Van Assche \cite{gerass} and the fact that there are only finitely many
eigenvalues in each gap was first proven in Cojuhari \cite{coj} and later rediscovered in
Teschl \cite{tosc}. Section~\ref{secTO} introduces the transformation operator and
proves the crucial decay estimate on its coefficients. This was first done by  Boutet de
Monvel and Egorova \cite{bdme} in the periodic case under the additional assumption
that all spectral gaps are open. We fix a problem in the original proof and at the same
time simplify and streamline the argument. Section~\ref{secSM} investigates the
scattering matrix. Our main result here is the reconstruction of the transmission coefficient
from the reflection coefficient, which was not known previously, even in the periodic case.
Section~\ref{secGLM} derives the Gel'fand-Levitan-Marchenko equation
and proves positivity of the Gel'fand-Levitan-Marchenko operator. In addition, we formulate
necessary conditions for the scattering data to uniquely determine our Jacobi operator.
Our final Section~\ref{secINV} shows that our necessary conditions for the scattering
data are also sufficient. It should be mentioned that, due to the lack of continuity
with respect to the spacial variable $n$, a significant change in the
strategy of the original proof in the continuous case from \cite{mar} is needed.

Our approach uses heavily the fact that the Baker-Akhiezer function is a meromorphic function
on the Riemann surface associated with the problem. This strategy gives a more streamlined
treatment and more elegant proofs even in the special cases which were previously
known. In this respect it is important to emphasize that, in contradistinction to the constant
background case, the upper sheet of our Riemann surface is not simply connected and
in particular not isomorphic to the unit disc.

\section{Quasi-periodic finite-gap operators and Riemann surfaces}
\label{secQP}

To set the stage let $\M$ be the Riemann surface associated with the following function
\begin{equation}
\Rg{z}, \qquad R_{2g+2}(z) = \prod_{j=0}^{2g+1} (z-E_j), \qquad
E_0 < E_1 < \cdots < E_{2g+1},
\end{equation}
$g\in \N$. $\M$ is a compact, hyperelliptic Riemann surface of genus $g$.
We will choose $\Rg{z}$ as the fixed branch
\begin{equation}
\Rg{z} = -\prod_{j=0}^{2g+1} \sqrt{z-E_j},
\end{equation}
where $\sqrt{.}$ is the standard root with branch cut along $(-\infty,0)$.

A point on $\M$ is denoted by 
$p = (z, \pm \Rg{z}) = (z, \pm)$, $z \in \C$, or $p = \infty_{\pm}$, and
the projection onto $\C \cup \{\infty\}$ by $\pi(p) = z$. 
The points $\{(E_{j}, 0), 0 \leq j \leq 2 g+1\} \subseteq \M$ are 
called branch points and the sets 
\begin{equation}
\Pi_{\pm} = \{ (z, \pm \Rg{z}) \mid z \in \C\backslash
\bigcup_{j=0}^g[E_{2j}, E_{2j+1}]\} \subset \M
\end{equation}
are called upper, lower sheet, respectively.

Let $\{a_j, b_j\}_{j=1}^g$ be loops on the surface $\M$ representing the
canonical generators of the fundamental group $\pi_1(\M)$. We require
$a_j$ to surround the points $E_{2j-1}$, $E_{2j}$ (thereby changing sheets
twice) and $b_j$ to surround $E_0$, $E_{2j-1}$ counter-clock wise on the
upper sheet, with pairwise intersection indices given by
\begin{equation}
a_i \circ a_j= b_i \circ b_j = 0, \qquad a_i \circ b_j = \delta_{ij},
\qquad 1 \leq i, j \leq g.
\end{equation}
The corresponding canonical basis $\{\zeta_j\}_{j=1}^g$ for the space of
holomorphic differentials can be constructed by
\begin{equation}
\underline{\zeta} = \sum_{j=1}^g \underline{c}(j)  
\frac{\pi^{j-1}d\pi}{R_{2g+2}^{1/2}},
\end{equation}
where the constants $\underline{c}(.)$ are given by
\[
c_j(k) = C_{jk}^{-1}, \qquad 
C_{jk} = \int_{a_k} \frac{\pi^{j-1}d\pi}{R_{2g+2}^{1/2}} =
2 \int_{E_{2k-1}}^{E_{2k}} \frac{z^{j-1}dz}{\Rg{z}} \in
\R.
\]
The differentials fulfill
\begin{equation}
\int_{a_j} \zeta_k = \delta_{j,k}, \qquad \int_{b_j} \zeta_k = \tau_{j,k}, 
\qquad \tau_{j,k} = \tau_{k, j}, \qquad 1 \leq j, k \leq g.
\end{equation}

Now pick $g$ numbers (the Dirichlet eigenvalues)
\be
(\hat{\mu}_j)_{j=1}^g = (\mu_j, \sigma_j)_{j=1}^g
\ee
whose projections lie in the spectral gaps, that is, $\mu_j\in[E_{2j-1},E_{2j}]$.
Associated with these numbers is the divisor $\dimuz$ which
is one at the points $\hat{\mu}_j$  and zero else. Using this divisor we
introduce
\begin{equation}
\ulz(p,n) = \hAmap(p) - \amap(\dimuz) - n\ul{A}_{\infty_-}(\infty_+) - \hvrc \in \C^g, \quad
\ulz(n) = \ulz(\infty_+,n),
\end{equation}
where $\vrc$ is the vector of Riemann constants
\begin{equation}
\hat{\Xi}_{p_0,j} = \frac{1- \sum_{k=1}^g \tau_{j,k}}{2},
\quad p_0=(E_0,0),
\end{equation}
and $\Amap$ ($\amap$) is Abel's map (for divisors). The hat indicates that we
regard it as a (single-valued) map from $\hat{M}$ (the fundamental polygon
associated with $\M$) to $\C^g$.
We recall that the function $\theta(\ulz(p,n))$ has precisely $g$ zeros
$\hmu_j(n)$ (with $\hmu_j(0)=\hmu_j$), where $\theta(\ul{z})$ is the
Riemann theta function of $\M$.

Then our Jacobi operator $H_q$ is given by
\bea \nn
a(n)^2 &=& \ti{a}^2 \frac{\theta(\ulz(n+1)) \theta(\ulz(n-1))}{\theta(
\ulz(n))^2},\\ \label{imfab}
b(n) &=& \tilde{b} + \sum_{j=1}^g c_j(g)
\frac{\partial}{\partial w_j} \ln\Big(\frac{\theta(\ul{w} +
\ulz(n)) }{\theta(\ul{w} + \ulz(n-1))}\Big) \Big|_{\ul{w}=0}.
\eea
The constants $\ti{a}$, $\tilde{b}$ depend only on the Riemann surface and
will be defined in the next section.

It is well known that the spectrum of $H_q$ is purely absolutely continuous and
consists of $g+1$ bands
\begin{equation}
\sig(H_q) = \bigcup_{j=0}^g [E_{2j},E_{2j+1}].
\end{equation}
For further information and proofs we refer to \cite{tjac}, Section~9.

\section{The Baker-Akhiezer function and the quasi-momentum map}
\label{secBA}

The Baker-Akhiezer function $\psi_q(p,n)=\psi_q(p,n,0)$ is given by
\begin{equation}
\psi_q(p,n,n_0) = \sqrt{\frac{\theta(\ulz(n_0-1))\theta(\ulz(n_0))}{\theta(\ulz(n-1))\theta(\ulz(n))}}
\frac{\theta(\ulz(p,n))}{\theta(\ulz(p,n_0))}
\: \exp \Big( (n-n_0) \int_{p_0}^p \hat{\omega}_{\infty_+,\infty_-} \Big),
\end{equation}
where $\omega_{\infty_+,\infty_-}$ is the normalized Abelian
differential of the third kind with simple poles at $\infty_\pm$ and residues
$\pm 1$, respectively. They are normalized such that $\psi_q(p,n_0,n_0)=1$.

The two branches
\begin{equation} \label{psi phi}
\psi_{q,\pm}(z,n) = \prod_{j=0}^{n-1} \phi_{q,\pm}(z,j),
\end{equation}
where (\cite{tjac}, (8.87)) 
\begin{equation} \label{def phi}
\phi_{q,\pm}(z,n) = \frac{1}{2a_q(n)}\biggl(z-b_q(n)+\sum_{j=1}^g
\frac{\hat R_j(n)}{z - \mu_j(n)} \pm \frac{R^{1/2}_{2g+2}(z)}{\prod_{j=1}^g (z - \mu_j(n))} \biggr), 
\end{equation}
\[
R_j(n) = 
\frac{R_{2g+1}^{1/2}(\mu_j(n))}{\prod_{k\ne j}(\mu_j(n)-\mu_k(n))}, \qquad
\hat R_j(n) = \sigma_j(n) R_j(n),
\]
of the Baker-Akhiezer function are
solutions of $\tau_q u = zu$, $z\in\C$, where $\tau_q$ is the difference
expression associated with $H_q$. However, the Wronskian
\be
W(\psi_{q,-}(z), \psi_{q,+}(z)) = \frac{R^{1/2}_{2g+2}(z)}{\prod_{j=1}^g (z-\mu_j)}
\ee
($\mu_j=\mu_j(0)$) shows that they are linearly dependent at the band edge
 $E_j$, $0 \leq j \leq 2g+1$.

The branch $\psi_{q,\sig_j}(z,n)$ has a first order pole at $\mu_j$ if $\mu_j$ is
away from the band edges
\be
\lim_{z\to\mu_j} (z-\mu_j) \psi_{q,\sig_j}(z,n) = 
\psi_{q,\sig_j}(\mu_j,n,1) \frac{\hat R_j(0)}{a_q(0)}
\ee
(use (\ref{def phi}) and $\psi_{q,\pm}(z,n)=\psi_{q,\pm}(z,n,1)\phi_{q,\pm}(z,0)$) and
both branches have a square root singularity if $\mu_j$ coincides with a
band edge $E_l$
\be
\lim_{z\to\mu_j} \sqrt{z-\mu_j} \psi_{q,\pm}(z,n) = \pm
\frac{\I^l \prod_{k\ne l}\sqrt{|E_l-E_k|}}{2a_q(0) \prod_{k\ne j}\sqrt{E_l-\mu_k}}
\psi_{q,+}(E_l,n,1)
\ee

\begin{lemma}     \label{le: char sol}
The solutions of $\tau_q u = zu$ can be characterized as follows.

(i) If $R_{2g+2}(z) \neq 0$, there exist two solutions satisfying
\be
\psi_{q,\pm}(z, n) = \theta_{\pm}(z,n) w(z)^{\pm n}, 
\quad w(z) = \exp \Big( \int_{p_0}^{(z,+)} \hat{\omega}_{\infty_+,\infty_-}
\Big),
\ee
with $\theta_{\pm}(z,n)$ quasi-periodic.    

(ii) If $R_{2g+2}(z) = 0$, $z=E_l$, there are two solutions satisfying
\be
\psi_q(E_l, n)= \psi_{q,+}(E_l,n) = \psi_{q,-}(E_l, n), \qquad 
\hat \psi_q(E_l, n) = \psi_q(E_l, n)(\hat \theta_l(n) + n),
\ee
where $\hat \theta_l(n)$ is quasi-periodic.
\end{lemma}

\begin{proof}
(ii). 
We construct a second linearly independent solution at $z = E = E_l$ using
(see \cite{tjac}, (1.50))
\be
s_q(E, n) = \lim_{z \rightarrow E} a_q(0) 
\frac{\psi_{q,+}(z,n) - \psi_{q,-}(z,n)}
{W(\psi_{q,-}(z), \psi_{q,+}(z))}, 
\ee
where $s_q(z, n)$ denotes the fundamental solution of $\tau_q u = z u$
with initial conditions $s_q(z,0)=0$, $s_q(z,1)=1$.
W.l.o.g.\ we assume that $E_l$ does not coincide with one of the Dirichlet
eigenvalues $\mu_j$ (otherwise shift the base point).
To derive an expression for $\psi_{q, \pm}(z)$ at $z = E + \epsilon^2$ we start with
\[
R^{1/2}_{2g+2}(z) = \epsilon (\tilde R + O(\epsilon^2)), \qquad
\tilde R = - \prod_{j \neq l}\sqrt{E - E_j}.
\]
Moreover,     
\[
W(\psi_{q,-}(z), \psi_{q,+}(z)) = \frac{\tilde R}{\prod_{j=1}^g (E-\mu_j)} \epsilon (1 + O(\epsilon^2))
\]
and for $p = (E + \epsilon^2, \pm)$ (see (\ref{ompminf}) below),
\[
\int_{p_0}^p \hat \omega_{\infty_+, \infty_-} =
\int_{p_0}^E \hat \omega 
\pm \beta \epsilon + O(\epsilon^3), \qquad
\beta = \frac{2 \prod_{j=1}^g(E - \lambda_j)}{\tilde R},
\]
\[
\underline{z}(p, n) = \underline{z}(E, n) \pm \underline{\gamma}\, \epsilon +
O(\epsilon^3), \qquad
\underline{\gamma} = \sum_{j=1}^g \underline{c}(j)
\frac{2E^{j-1}}{\tilde R},
\]
and 
\[
\theta (\underline{z}(p, n)) = \theta(\underline{z}(E, n)) \pm
\frac{\partial \theta}{\partial \underline{z}} (\underline{z}(E, n))\,
\underline{\gamma}\, \epsilon + O(\epsilon^3).
\]
Using this to evaluate the limit $\eps\to 0$ shows
\[
s_q(E,n) = 2a_q(0) \prod_{j=1}^g \frac{E-\mu_j}{E-\lam_j} \hat \psi_q(E, n) =
\psi_q(E, n) (\hat \theta(n) + n),
\]
where
\[
\hat \theta (n) = \frac{1}{\prod_{j=1}^g (E - \lambda_j)}
\sum_{j,k=1}^g E^j c_k(j) 
\frac{\partial}{\partial w_k} \ln \theta(\ul{z}(E,n) + \ul{w}),
\]
and finishes the proof.
\end{proof}

\begin{remark}
(i). Since $\psi_q(z,n)$ has a singularity if $z=\mu_j$ the solutions in
Lemma~\ref{le: char sol} are not well-defined for those $z$. However,
you can either remove the singularities of $\psi_q(z,n)$ or choose
a different normalization point $n_0\ne 0$ to see that solutions of
the above type exist for every $z$.\\
(ii). In the periodic case Floquet theory tells you that there are two
possible cases at a band edge: Either two (linearly independent)
periodic solutions or one periodic and one linearly growing solution.
The above lemma shows that the first case happens if the corresponding
gap is closed and the second if the gap is open.
\end{remark}

To understand the properties of $\psi_{q,\pm}(z,n)$ we need to investigate the quasi-momentum map
\begin{equation}
w(z) = \exp \Big( \int_{p_0}^p \hat{\omega}_{\infty_+,\infty_-} \Big),
\quad p=(z,+).
\end{equation}
The differential $\omega_{\infty_+,\infty_-}$ is given by
\begin{equation} \label{ompminf}
\om_{\infty_+,\infty_-} = \frac{\prod_{j=1}^g (\pi - \lam_j)}{R_{2g+2}^{1/2}} d\pi,
\end{equation}
where the constants $\lam_j$ have to be determined from the normalization
\begin{equation} \label{ominfpmpl}
\int_{a_j} \om_{\infty_+,\infty_-} = 2 \int\limits_{E_{2j-1}}^{E_{2j}} 
\frac{\prod_{j=1}^g (z - \lam_j)}{\Rg{z}} dz = 0,
\end{equation}
which shows $\lam_j \in (E_{2j-1},E_{2j})$.

Since $\lambda_j \in (E_{2j-1}, E_{2j})$ the integrand
is a Herglotz function and admits the following representation (c.f.\
\cite{tjac}, Appendix B)
\be
 \frac{\prod_{j=1}^g (z - \lambda_j)} {\Rg{z}}
 = \int_{- \infty}^{\infty}\frac{1}{\lambda - z}
d \tilde \mu(\lambda)
\ee
with the probability measure
\be 
d\tilde\mu(\lambda) = \frac{\prod_{j=1}^g(\lambda - \lambda_j)}
{\pi\I\Rg{\lam}}
\chi_{\sigma(H_q)}(\lambda) d\lambda.
\ee
Hence
\bea             \nonumber
g(z,\infty) &=& \int_{p_0}^p \om_{\infty_+,\infty_-} = \int_{E_0}^z 
\int_{- \infty}^{\infty}\frac{1}{\lambda - \zeta}
d \tilde \mu(\lambda) d\zeta     \\
&=& \int_{- \infty}^{\infty} 
\ln\left(\frac{\lambda - E_0}{\lambda - z} \right) d \tilde \mu(\lambda).
\eea
In particular, note that $-\mbox{Re}(g(z,\infty))$ is the Green's function of the
upper sheet $\Pi_+$ with pole at $\infty_+$ and $\tilde{\mu}$ is the
equilibrium measure of the spectrum (see \cite{tsu}, Thm.~III.37).
We will abbreviate $g(z)=g(z,\infty)$.

The asymptotic expansion of $\exp(g(z))$ is given by (\cite{tjac}, (9.42))
\be     \label{asymp w}
\exp \bigg(\int_{p_0}^p \hat \omega_{\infty_+, \infty_-}\bigg)
= -  \frac{\ti{a}}{z}  \Big(1 + \frac{\ti{b}}{z} + O(\frac{1}{z^2})\Big),
\qquad z \rightarrow \infty,
\ee
where $\ti{a}$ is the capacity of the spectrum and
\be
\ti{b} = \frac{1}{2} \sum_{j=0}^{2g+2} E_j - \sum_{j=1}^g \lambda_j.
\ee

\begin{theorem}  
The map $g$ is a bijection from the upper (resp. lower) half plane 
$\C^{\pm}=\{z \in \C \mid \pm \im(z)>0\}$ to
\be
S^{\pm} = \{z \in \C \mid \pm \re(z) < 0, 0 < \im(z) < \pi \}
\backslash \bigcup_{j=1}^g [g(\lambda_j), g(E_{2j+1})]
\ee
such that $\sigma(H_q) = \{z \mid \re(z) = 0\}$.
\end{theorem}

\begin{proof}
By the Herglotz property of its integrand, the function $g(z,\infty)$
satisfies the conditions of \cite{par}, Theorem 1(b) in Chapter VI,
which shows that it is one-to-one.

To prove that $g(z,\infty)$ is surjective, it suffices to show that the boundary of 
$\C^+$ is mapped to the boundary of $S^+$. Note that
$g(\lambda)$ is negative for $\lambda < E_0$ and purely imaginary for 
$\lambda \in [E_0, E_1]$. At $E_1$, the real part starts to decrease from 
zero until it hits its minimum at $\lambda_1$ and increases again until it
becomes 0 at $E_2$ (since all $a$-periods are zero), while the imaginary part remains constant.
Proceeding like this we move along the boundary of $S^+$ as $\lambda$ moves
along the real line. For $\lambda > E_{2g+1}$, $g(\lambda)$ is again
negative.
\end{proof}

\begin{remark}
In the special case where $H_q$ is periodic the quasi-momentum is given by 
$w(z)=\exp(\I N^{-1}\arccos \Delta(z))$, where $\Delta(z)$ is the Floquet 
discriminant, and our result is due to \cite{perco}.
\end{remark}

Therefore the map  
\bea  \nonumber
w: \C^{\pm} &\to& 
W^{\pm}=\{w \in \C \mid |w|< 1, \pm \mbox{Im}(w) > 0\} \backslash 
\bigcup_{j=1}^g [w(\lambda_j), w(E_{2j+1})]           \\
z &\mapsto& \exp(g(z))
\eea
is bijective. Denote $W = W^+ \cup W^- \cup (-1, 1)$, $W_0= W \backslash\{0\}$.
If we identify corresponding points on the slits $[w(\lambda_j), w(E_{2j+1})]$ we
obtain a Riemann surface $\W$ which is isomorphic to the upper sheet $\Pi_+$.

\begin{remark}
In \cite{perco} the largest band edge $E_{2g+1}$ is chosen for $p_0$ and $w$ will
map $\C^\pm \to W^\mp$ in this case. Moreover, in the periodic case the
slits $[w(\lambda_j), w(E_{2j+1})]$ appear at equal angles $\frac{2\pi}{N}$, where
$N$ is the period.
\end{remark}

Since $z \mapsto w(z)= \exp(g(z))$ is a bijection, we consider the
functions $\psi_{q,\pm}$ as functions of the new parameter
$w$ whenever convenient. For notational simplicity we will write
$\psi_{q,\pm}(w,n)$ for $\psi_{q,\pm}(\lambda(w),n)$ and similarly
for other quantities.
The functions $\psi_{q,\pm}(w,n)$ are meromorphic in $\W$
and continuous up to the boundary with the only possible singularities
at the images of the Dirichlet eigenvalues $w(\mu_j)$ and at $0$.
More precisely, denote by $M_{\pm}$ the sets of poles (and square root 
singularities if $\mu_j = E_l$) of the Weyl $m$-functions 
$\tilde m_{\pm}(\lambda)$, i.e.\ $M_+ \cup M_- = \{\mu_j\}_{j=1}^g$
(see (\ref{psi phi}) and \cite{tjac}, Section~2.1). Note that
$\mu_j \in M_+ \cap M_-$ if and only if $\mu_j = E_l$. Then

\begin{itemize}
\item[(B1)] $\psi_{q,\pm}(w,n)$ are holomorphic in
$\W \backslash (\{ w(\mu_j) \}_{j=1}^g \cup \{ 0\})$ and 
continuous on $\partial W\backslash\{w(\mu_j)\}$.
\item[(B2)] $\psi_{q,\pm}(w,n)$ has a simple pole at $w(\mu_j)$ if
$\mu_j\in M_\pm \backslash \{E_l\}$, no pole if $\mu_j\not\in M_\pm$,
and if $\mu_j = E_l$,
$$
\psi_{q,\pm}(w, n) = \pm 
\frac{\I^l C(n)}{w-w_l} + O(1),
$$
where $C(n)$ is bounded and real.
\item[(B3)] $\overline{\psi_{q,\pm}(w, n)}= \psi_{q,\mp}(w, n)$ for $|w|=1$. 
\item[(B4)]  At $w=0$ the following asymptotics hold 
$$
\psi_{q,\pm}(w,n) =  (-1)^n  
\left(\frac{\prod_{m=0}^{n-1}a_q(m)}{\tilde a^{n}}\right)^{\pm 1} w^{\pm n} (1+O(w)).
$$
\end{itemize}

By Section~2.5 of \cite{tjac} the vector valued functions 
\be
\underline{U}(\lambda, n) =  
\sqrt{\frac{1}{4a_q(0)^2\pi \mbox{\rm Im}(\tilde m_+(\lambda))}}
\left( \begin{array}{c} 
\psi_{q,+}(\lambda, n)  \\
\psi_{q,-}(\lambda, n)  \\
\end{array} \right)
\ee
form an orthonormal basis for the Hilbert space 
$L^2(\sigma(H_q), \C^2, d\lambda)$. The Weyl $m$-functions $\ti{m}_\pm(z)$  
satisfy (see \cite{tjac}, eq. (8.27))
\be
\im(\ti{m}_\pm(\lam)) = \frac{\mp\Rg{\lam}}{2 \I a_q(0)^2 \prod_{j=1}^g (\lam - \mu_j)},
\qquad \lam\in\sig(H_q).
\ee

Using our map $w(z) =\exp(\int_{p_0}^{(z,+)} \hat \omega_{\infty_+, \infty_-})$ we can
transform this into an orthogonal basis on the unit circle.

\begin{lemma}                    \label{lem orth}
Both functions $\psi_{q,+}(w, n)$ and $\psi_{q,-}(w, n)$ form
orthonormal bases in the Hilbert space $L^2(S^1, \frac{1}{2\pi\I} d\omega)$,
where
\begin{equation}
d\omega(w) = \prod_{j=1}^g \frac{\lambda(w)- \mu_j}
{\lambda(w)-\lambda_j} \frac{dw}{w}.
\end{equation}
\end{lemma}

\begin{proof}
Just use
\be
\frac{d w}{dz} = w \frac{\prod_{j=1}^g(z - \lambda_j)}{\Rg{z}}.
\ee
\end{proof}

Observe that $d\omega$ is meromorphic on $\W$ with a simple pole at $w=0$. In particular,
there are no poles at $w(\lam_j)$.

\begin{remark}
In the periodic case we have
\be                                      \label{norm m_+ g}
\| \psi_{p, \pm}(\lambda)\|_{N}^2 := 
\sum_{n=1}^{N}|\psi_{p, \pm}(\lambda, n)|^2
= N \prod_{j=1}^{N-1} \frac{\lambda-\lambda_j}{\lambda- \mu_j}.
\ee
\end{remark}

\section{Existence of Jost solutions}
\label{secJS}

After we have these preparations out of our way, we come to
the study of short-range perturbations $H$ of $H_q$ associated with sequences $a$, $b$
satisfying $a(n) \rightarrow a_q(n)$ and $b(n) \rightarrow b_q(n)$ as $|n|
\rightarrow \infty$.  More precisely, we will make the following assumption throughout
this paper.

\begin{hypo}            \label{hypo:quasi}
Let $H$ be a perturbation such that
\begin{equation}                         \label{hypo}
\sum_{n \in \mathbb{Z}} |n| \Big(|a(n) - a_q(n)| + |b(n) - b_q(n)| \Big) <
\infty.
\end{equation}
\end{hypo}

We first establish existence of Jost solutions, that is solutions of the
perturbed operator which asymptotically look like the Baker-Akhiezer solutions.

\begin{theorem} \label{thmjost}
Assume (H.\ref{hypo:quasi}). Then there exist solutions
$\psi_{\pm}(z, .)$, $z \in \C$, of $\tau \psi = z \psi$ satisfying
\begin{equation}                         \label{perturbed sol}
\lim_{n \rightarrow \pm \infty}
|w(z)^{\mp n} (\psi_{\pm}(z, n) - \psi_{q,\pm}(z, n))| = 0,
\end{equation}
where $\psi_{q,\pm}(z, .)$ are the Baker-Akhiezer functions.
Moreover, $\psi_{\pm}(z, .)$ are continuous 
(resp.\ holomorphic) with respect to $z$ whenever $\psi_{q,\pm}(z, .)$ are and
inherit the properties $(B1)$ and $(B2)$, where now $\psi_\pm(z, n) = \frac{\I^l C_\pm(n)}{\sqrt{z-\mu_j}} + O(1)$ and
$(B4)$ has to be replaced by
\begin{equation}                 \label{B4jost}
\psi_\pm(z,n) =   A_\pm(0) \left(\frac{\prod_{m=0}^{n-1}a(m)}{z^n}\right)^{\pm 1} 
(1 + (B_\pm(0) \pm \sum_{j=1}^n b(j- {\scriptstyle{0 \atop 1}}))\frac{1}{z} + O(\frac{1}{z^2})),
\end{equation}
where
\bea \nonumber
A_+(n) &=& \prod_{j=n}^{\infty} \frac{a_q(j)}{a(j)}, \qquad
B_+(n)= \sum_{m=n+1}^\infty (b_q(m)-b(m)), \\
A_-(n) &=& \prod_{j=- \infty}^{n-1} \frac{a_q(j)}{a(j)}, \qquad
B_-(n) = \sum_{m=-\infty}^{n-1} (b_q(m)-b(m)).
\eea
\end{theorem}

\begin{proof}
The proof can be done as in the periodic case (see e.g., \cite{gerass}, \cite{tosc} or
\cite{tjac}, Section 7.5). 
The only problem is to show that the second solution at a band edge grows at
most linearly. In the periodic case this follows from Floquet theory, here
we just use Lemma~\ref{le: char sol}.
\end{proof}

 
 From this result we obtain a complete characterization of the spectrum of $H$.

\begin{theorem} 
Assume (H.\ref{hypo:quasi}). Then we have $\sig_{ess}(H)=\sig(H_q)$, the
point spectrum of $H$ is finite and confined to the spectral gaps of
$H_q$, that is, $\sig_p(H) \subset \R\backslash\sig(H_q)$. Furthermore,
the essential spectrum of $H$ is purely absolutely continuous.
\end{theorem}

\begin{proof}
Again the proof can be done as in the periodic case
(see e.g., \cite{tosc} or \cite{tjac}, Section 7.5).
\end{proof}

\section{The transformation operator}
\label{secTO}

We define the kernel of the transformation operator as the Fourier
coefficients of the Jost solutions $\psi_{\pm}(w,n)$ 
with respect to the orthonormal system given in Lemma~\ref{lem orth}, 
$\{\psi_{q,\pm}(w,n)\}_{n\in \Z}$,  
\be                             \label{de:K qp}
K_{\pm}(n,m) := \frac{1}{2\pi \I} \int_{|w|=1} \psi_{\pm}(w,n) \psi_{q,\mp}(w, m)
d\omega(w).
\ee
By the Cauchy theorem, this integral equals the residue at $w=0$,
\be
K_{\pm}(n,m) = \res_{0} \frac{1}{w} \psi_{\pm}(w,n) \psi_{q,\mp}(w, m).
\ee
In particular, since $\psi_{\pm}(w,n) \psi_{q,\mp}(w, m) = O(w^{\pm(n-m)})$, we
conclude
\be
K_\pm(n,m) =0, \qquad \pm(m-n) < 0.
\ee

\begin{lemma}                     \label{le:jost kernel q}
Assume H.\ref{hypo:quasi}.
The Jost solutions $\psi_{\pm}(w,n)$ can be represented as  
\be                     \label{psi repr}
\psi_{\pm}(w,n) = \sum_{m=n}^{\pm \infty} K_{\pm}(n,m)
\psi_{q,\pm}(w, m), \qquad |w|=1,
\ee
where the kernels $K_{\pm}(n, .)$ satisfy 
$K_{\pm}(n,m) = 0$ for $\pm m < \pm n$ and  
\be                      \label{estimate K q}
|K_{\pm}(n,m)| \leq C \sum_{j=[\frac{m+n}{2}] \pm 1}^{\pm \infty}
\Big(|a(j)-a_q(j)| + |b(j)-b_q(j)|\Big), \quad \pm m > \pm n.
\ee
The constant $C$ depends only on $H_q$ and the value of the sum in 
(\ref{hypo}).
\end{lemma}

\begin{proof}  
We prove the estimate for $K_+(n,m)$ and omit "$+$" and 
"$z$" whenever possible. Define $\varphi(n) = \psi(n) K(n,n)^{-1}$,
then $\varphi$ fulfills 
\be
\varphi(n) = \psi_q(n) + \sum_{m=n+1}^{\infty} J(n,m) \varphi(m),
\ee
where
\be
J(z,n,m)= \tilde a(m-1) \frac{s_q(z,n,m-1)}{a_q(m-1)} 
+ \tilde b(m) \frac{s_q(z,n,m)}{a_q(m)}
\ee
with the abbreviation
\be                      
\tilde a(m) = \frac{a(m)^2}{a_q(m)} - a_q(m), \qquad 
\tilde b(m) = b(m) - b_q(m).
\ee
On the other hand, $\varphi(n)$ is given by 
\[
\varphi(n) = \sum_{m=n}^{\infty}\kappa(n,m) \psi_q(m), \qquad
\kappa(n,m) = \frac{K(n,m)}{K(n,n)},
\]
therefore 
\be                           \label{varphi 2}
\sum_{m=n}^{\infty}\kappa(n,m) \psi_q(m) = 
\sum_{m=n+1}^{\infty} J(n,m) \psi_q(m) 
+ \sum_{m=n+1}^{\infty} \sum_{l=m+1}^{\infty}J(n,m)
\kappa(m,l) \psi_q(l).    
\ee
Multiplying both sides of (\ref{varphi 2}) by $\psi_{q,-}(k)$ and
integrating over the unit circle yields
\be                              \label{die1}
\kappa(n,k) = \sum_{m=n+1}^{\infty} \Gamma(n,m,m,k) +
\sum_{m=n+1}^{\infty} \sum_{l=n+1}^{\infty} \Gamma(n,m,l,k) \kappa(m,l), 
\ee
where
\be
\Gamma(n,m,l,k) = \frac{1}{2 \pi \I} \int_{|w|=1} J(w,n,m)
\psi_{q,+}(w,l) \psi_{q,-}(w, k) d\omega(w).
\ee
Using \cite{tjac}, (1.50),
\be
\frac{s_q(n,m)}{a(m)} = \frac{\psi_{q, +}(m)\psi_{q, -}(n) - 
\psi_{q, +}(n)\psi_{q, -}(m)}{W(\psi_{q, +}, \psi_{q, -})},
\ee
we obtain  
\be
\Gamma(n,m,l,k) = \tilde b(m) \Gamma_q(n,m,l,k) + \tilde a(m)
\Gamma_q (n,m-1,l,k) 
\ee
with
\bea \nonumber                  \label{Gamma_0}
\Gamma_q(n,m,l,k) &=& \Gamma_0(m,n,l,k) - \Gamma_0 (n,m,l,k), \\ \nonumber
\Gamma_0(n,m,l,k) &=& \frac{1}{2 \pi \I} \int_{w(\gam)} 
\frac{\psi_{q,+}(w,n) \psi_{q, -}(w, m)
\psi_{q,+}(w,l) \psi_{q,-}(w, k)}{W(\psi_{q, +}(w), \psi_{q, -}(w))} 
d\omega(w)             \\ \nonumber
&=& \frac{1}{2 \pi \I} \int_{\gam} 
\frac{\psi_{q, +}(z, n)\psi_{q, -}(z, m)
\psi_{q,+}(z,l) \psi_{q,-}(z, k)}{W(\psi_{q, +}(z), \psi_{q, -}(z))} 
\frac{\prod(z - \mu_j)}{\Rg{z}} dz   \\  
&=& \frac{1}{2 \pi \I} \int_{\gam} 
\frac{\psi_{q, +}(z, n)\psi_{q, -}(z, m)
\psi_{q,+}(z,l) \psi_{q,-}(z, k)}{W(\psi_{q, +}(z), \psi_{q, -}(z))^2} dz.
\eea
Here $\gam$ is a path on the upper sheet encircling the spectrum.
The integrand of $\Gamma_0$ is meromorphic on the Riemann surface $\M$ with
poles of order one at $E_j$ and poles of order $O(z^{\pm(n-m+l-k)-2})$ near 
$\infty_{\pm}$ (there are no poles at the Dirichlet eigenvalues $\mu_j$ ).
We apply the residue theorem twice, first on the side of 
$\gamma$ including $\infty_+$, then on the other side including the
spectrum (and thus $\infty_-$)
\bea    \nn                                
\Gamma_0(n,m,l,k) &=& - \res_{\infty_+} \frac{\psi_{q, +}(n)
\psi_{q, -}(m) \psi_{q,+}(l) \psi_{q,-}(k)}
{W(\psi_{q, +}, \psi_{q, -})^2}  \\
&=& \Big(\res_{\infty_-} + \sum_{j=0}^{2g+1}\res_{E_j}\Big)
\left( \frac{\psi_{q, +}(n) \psi_{q, -}(m) \psi_{q,+}(l) \psi_{q,-}(k)}
{W(\psi_{q, +}, \psi_{q, -})^2} \right).  
\eea
The order of the poles at $\infty_{\pm}$ implies
\be                             \nonumber
\Gamma_0(n,m,l,k) = \left\{ 
\array{cc}      
\sum\limits_{j=0}^{2g+1}\res_{E_j} 
\frac{\psi_{q, +}(n) \psi_{q, -}(m) \psi_{q,+}(l) \psi_{q,-}(k)}
{W(\psi_{q, +}, \psi_{q, -})^2} & n-m+l-k < 0 \\
0                       &  n-m+l-k \geq 0, \\
\endarray
\right.          
\ee
which shows that $\Gamma_0(n,m,l,k)$ is real and bounded since 
$\psi_{q,+}(E, .)=\psi_{q,-}(E, .)$ are (if $\mu_j=E_l$, use (B2)).
Together with (\ref{Gamma_0}) this yields
\be
\Gamma_0(n,m,l,k) = - \overline{\Gamma_0(m,n,k,l)} = - \Gamma_0(m,n,k,l)
= - \Gamma_0(n,m,k,l).
\ee
Moreover,  
\bea  \nonumber
\Gamma_q(n,m,l,k) &=& 0, \qquad l-k \geq |m-n|, \\ 
\Gamma_q(n,m,l,k) &=& - \Gamma_q(m,n,k,l) \ =\ \Gamma_q(n,m,k,l),
\eea
which then implies
\be         \label{Gamma 3}                     
\Gamma_q(n,m,l,k) = \left\{ \!\!\!
\array{cc}      
\mbox{sign}(n-m) \!\!\sum\limits_{j=0}^{2g+1} \!\! \res_{E_j} \!
\frac{\psi_{q, +}(n) \psi_{q, -}(m) \psi_{q,+}(l) \psi_{q,-}(k)}
{W(\psi_{q, +}, \psi_{q, -})^2}\!\! & |l-k| < |m-n| \\
0  &  |l-k| \geq |m-n| \\
\endarray
\right.          
\ee
and $\Gamma(n,m,l,k)=0$ for $|l-k|\geq m-n$ if $m>n$. Note that
the residue at $E_j$ is given by
\be			\label{res Gamma}
\frac{2 \prod_{\ell=1}^g (E_j -\mu_\ell)^2}{\prod_{\ell\ne j} (E_j-E_\ell)}
\psi_q(E_j,n) \psi_q(E_j,m) \psi_q(E_j,l) \psi_q(E_j,k).
\ee
Now we obtain for 
$\kappa(n,k)$
\bea       \nn                  
\kappa(n,k) &=& 
\sum_{m=n+1}^{\infty} \Gamma(n,m,m,k) +
\sum_{m=n+1}^{\infty} \sum_{l=m+1}^{\infty} \Gamma(n,m,l,k)
\kappa(m,l) \\            \label{Gamma 2}    
&=& \sum_{m=[\frac{n+k}{2}]+1}^{\infty} \Gamma(n,m,m,k) +
\sum_{m=n+1}^{\infty} \sum_{l=n+k-m+1}^{m+k-n-1} \Gamma(n,m,l,k)
\kappa(m,l),
\eea
since $\Gamma(n,m,m,k) \neq 0$ only if $|m-k|<m-n$ implying
$m>\frac{n+k}{2}$. In the third sum of (\ref{Gamma 2}) we need that 
$|m + \delta - k| < m-n$ for $\delta \geq 1$ which yields 
$\delta < k-n$ and $\delta > n+k-2m$.
Two remarks might be in order: $m+k-n-1 \geq n+k-m+1$ since
$m-n\geq n-m+2$, and the starting point $l=n+k-m+1$ of the third sum actually 
has a lower limit, namely $m \leq \frac{n+k}{2}$, since we require
$l \geq m +1$ for $\kappa(m,l)\neq 0, 1$.
Note that
\bea                                                    \nonumber
\Bigg|\sum_{m=[\frac{n+k}{2}]+1}^{\infty} \Gamma(n,m,m,k) \Bigg|
&\leq& D \sum_{m=[\frac{n+k}{2}]+1}^{\infty} |\tilde b(m) + \tilde a(m)|
=: \hat q(\tfrac{n+k}{2}),  \\
\Bigg| \sum_{l=n+k-m+1}^{m+k-n-1} |\Gamma(n,m,l,k)| \Bigg|  &\leq&
D\, (m-n-1) |\tilde b(m) + \tilde a(m)| =: \hat c(m)  \in \ell^1(\Z),    
\nonumber
\eea
where $D$ is the estimate provided by (\ref{Gamma 3}), (\ref{res Gamma}).
We set up the following iteration procedure 
\bea             \nonumber
\kappa_0(n,k) &=& \sum_{m=[\frac{n+k}{2}]+1}^{\infty} \Gamma(n,m,m,k), 
\\   
\kappa_j(n,k) &=& \sum_{m=n+1}^{\infty} 
\sum_{l=n+k-m+1}^{m+k-n-1} \Gamma(n,m,l,k) \kappa_{j-1}(m,l).
\eea
Then using induction one has
\be
| \kappa_j(n,k)| \leq \hat q(\tfrac{n+k}{2}) 
\frac{\left( \sum_{m=n+1}^{\infty} \hat c(m) \right)^j}{j!}
\ee
and hence the iteration converges and implies the estimate
\be                              \label{kappa 2}
|\kappa(n,k)| = \Bigg|\sum_{j=0}^{\infty}\kappa_j(n,k) \Bigg| \leq
\hat q(\tfrac{n+k}{2}) 
\exp \Bigg(\sum_{m=n+1}^{\infty} \hat c(m) \Bigg). 
\ee 
\end{proof}

Associated with $K_{\pm}(n,m)$ is the operator 
\be
(\mathcal{K_{\pm}} f)(n) = \sum_{m=n}^{\pm \infty} K_{\pm} (n,m) f(m), 
\qquad  f \in \ell_{\pm}^{\infty}(\Z, \C),
\ee
which acts as a transformation operator for
the pair $\tau$, $\tau_q$. 

\begin{theorem} 
Let $\tau_q$ and $\tau$ be the quasi-periodic and perturbed 
Jacobi difference expression, respectively. Then \label{th:transf op q}
\be
\tau \mathcal{K_{\pm}}f = \mathcal{K_{\pm}} \tau_q f, \qquad f \in
\ell_{\pm}^{\infty}(\Z, \C).
\ee
\end{theorem}

\begin{proof}
It suffices to show that $H K_{\pm} = K_{\pm} H_q$.
\bea         \nonumber
H K_{\pm}(n,m) &=& 
\frac{1}{2\pi \I} \int_{|w|=1} H \psi_{\pm}(w,n) 
\psi_{q,\mp}(w, m) d\omega(w)\\         \nonumber
&=& \frac{1}{2\pi \I} \int_{|w|=1} \lambda(w) \psi_{\pm}(w,n) 
\psi_{q,\mp}(w, m)d\omega(w)\\           
&=& \frac{1}{2\pi \I} \int_{|w|=1} \psi_{\pm}(w,n) 
H_q \psi_{q,\mp}(w, m) d\omega(w).      \label{H K = K H_q}
\eea
\end{proof}

\begin{lemma}      \label{coro K+K- q}
For $n \in \Z$ we have
\bea                                             \label{a a_q}
\frac{a(n)}{a_q(n)} &=& \frac{K_+(n+1, n+1)}{K_+(n,n)} 
\ =\  \frac{K_-(n,n)}{K_-(n+1, n+1)},  \\          \nonumber
b(n) - b_q(n) &=& a_q(n)\frac{K_+(n, n+1)}{K_+(n,n)}
-a_q(n-1)\frac{K_+(n-1, n)}{K_+(n-1,n-1)}  \\           \nonumber
&=& a_q(n-1)\frac{K_-(n, n-1)}{K_-(n,n)}
-a_q(n)\frac{K_-(n+1, n)}{K_-(n+1,n+1)}, 
\eea
\end{lemma}

\begin{proof}
Consider the equation of the transformation
operator $H \mathcal{K_{\pm}} = \mathcal{K_{\pm}} H_q$, 
which is equivalent to (c.f.\ (\ref{H K = K H_q}))
\[
a(n-1)K_{\pm}(n-1,m) + b(n)K_{\pm}(n,m) + a(n)K_{\pm}(n+1,m) =
\]
\[
= a_q(m-1)K_{\pm}(n,m-1) + b_q(m)K_{\pm}(n,m) + a_q(m)K_{\pm}(n,m+1).
\]
Evaluating at $m=n$ we obtain the first equation and at $m=n\mp 1$
the second.
\end{proof}

In particular, observe
\be                                             \label{K(n,n)}
K_\pm(n,n) = A_\pm(n), \qquad
K_\pm(n,n\pm 1) = \frac{A_\pm(n)}{a_q(n- {0 \atop 1})} B_\pm(n).
\ee

\section{The scattering matrix}
\label{secSM}

Let $H_q$ be a given quasi-periodic Jacobi operator and $H$ a perturbation 
of $H_q$ satisfying Hypothesis H.\ref{hypo:quasi}. 
To set up scattering theory for the pair $(H, H_q)$ we proceed
as usual.

The Wronskian of our Jost functions can be evaluated as $n\to\pm\infty$ and
is given by
\be             \label{W q}
W(\psi_{\pm}(\lambda), \overline{\psi_{\pm}(\lambda)}) = 
W_q(\psi_{q,\pm}(\lambda), \psi_{q,\mp}(\lambda)) = 
\mp \frac {\Rg{\lambda}}
{\prod_{j=1}^g(\lambda - \mu_j)}, \qquad \lambda \in \sigma(H_q).
\ee
Hence $\psi_{\pm}(\lambda)$, 
$\overline{\psi_{\pm}(\lambda)}$ are linearly independent for
$\lambda$ in the interior of $\sigma(H_q)$ and we consider the scattering
relations
\be                      \label{scat rel q}
\psi_{\pm}(\lambda,n) = 
\alpha(\lambda) \overline{\psi_{\mp}(\lambda,n)} 
+ \beta_{\mp}(\lambda)\psi_{\mp}(\lambda,n),
\qquad \lambda \in \sigma(H_q),
\ee
where
\bea                   \label{alpha quasi}                   
\alpha(\lambda) &=& \frac{W(\psi_{\mp}(\lambda), 
\psi_{\pm}(\lambda))}
{W(\psi_{\mp}(\lambda), \overline{\psi_{\mp}(\lambda)})}\ = \
\frac {\prod_{j=1}^g(\lambda - \mu_j)}
{\Rg{\lambda}}
W(\psi_-(\lambda), \psi_+(\lambda)), \\
\beta_{\pm}(\lambda) &=& \frac{W(\psi_{\mp}(\lambda), 
\overline{\psi_{\pm}(\lambda)})}
{W(\psi_{\pm}(\lambda), \overline{\psi_{\pm}(\lambda)})}\ 
= \ \mp \frac {\prod_{j=1}^g(\lambda - \mu_j)}{\Rg{\lambda}}      \nonumber
W(\psi_{\mp}(\lambda), \overline{\psi_{\pm}(\lambda)}). 
\eea
While $\alpha(\lambda)$ is only defined for $\lambda \in \sigma(H_q)$, 
(\ref{alpha quasi}) may be used as a definition for 
$\lambda \in \C\backslash \{E_j\}$. 
Therefore $\alpha(w)$ can be continued as a
holomorphic function on $\W$ and it is continuous up to the
boundary except possibly at the band edges.

\begin{remark}
Note that $\alpha(\lambda)$ does not depend on the normalization of 
$\psi_{\pm}(\lambda)$ at the base point $n_0=0$ whereas
$\beta_{\pm}=\beta_{\pm,0}$ does. Using 
$\psi_{\pm}(z,n,n_0)= \psi_{q,\pm}(z,n_0)^{-1}\psi_{\pm}(z,n)$ and
\[
W((\psi_+(\lambda), \psi_-(\lambda)) = \prod_{j=1}^g \frac{\lambda - \mu_j(n_0)}
{\lambda - \mu_j} W((\psi_+(\lambda, .,n_0), \psi_-(\lambda, ., n_0))
\]
we see
\begin{equation}
\beta_{\pm,0}(\lambda) = \frac{\psi_{q, \mp}(\lambda, n_0)}
{\psi_{q, \pm}(\lambda, n_0)} \beta_{\pm,n_0}(\lambda).
\end{equation}
\end{remark}
 
A direct calculation shows  
\be                      \label{alpha beta 1}
\alpha(\overline w) = \overline{\alpha(w)}, \quad
\beta_{\pm}(\overline w) = \overline{\beta_{\pm}(w)} = - \beta_{\mp}(w)
\ee
and the Pl\"ucker identity (c.f.\ \cite{tjac}, (2.169)) implies 
\be                      \label{alpha beta 2}
|\alpha(w)|^2 = 1 + |\beta_{\pm}(w)|^2, \qquad |w|=1.  
\ee

We will denote the eigenvalues of $H$ by 
\be
\sigma_p(H) = \{ \rho_j \}_{j=1}^q.
\ee
Our next aim is to study the behavior of $\alpha(\lambda)$ at the
eigenvalues $\rho_j$, therefore we modify the Jost solutions 
$\psi_{\pm}(\lambda, n)$ according to their poles at $\mu_j$ and define
the following eigenfunctions $\hat{\psi}_{\pm}(\lambda, .)$ 
\be    \label{def eigen}
\hat{\psi}_+(\lambda, .) = \prod_{\mu_l \in M_+}
(\lambda - \mu_l) \,\psi_+(\lambda, .), \quad
\hat{\psi}_-(\lambda, .) = \prod_{\mu_l \in M_- \backslash \{E_j\}}
(\lambda - \mu_l) \,\psi_-(\lambda, .).
\ee
Define $\hat{\psi}_{q, \pm}(\lambda, .)$ accordingly. 
Moreover, $\hat{\psi}_{\pm}(\rho_j, n) = c_j^{\pm} \hat{\psi}_{\mp}(\rho_j, n)$ with $c_j^+ c_j^- = 1$. 
The {\bf norming constants} $\gamma_{\pm,j}$ are defined
by
\be
\frac{1}{\gamma_{\pm,j}} := \sum_{m \in \Z} |\hat{\psi}_{\pm}(\rho_j, m)|^2.
\ee

To compute the derivative of $\alpha(\lambda)$ at $\rho_j$, note that
\be
\alpha(\lambda)   
= \frac{W(\hat{\psi}_{-}(\lambda), \hat{\psi}_{+}(\lambda))}
{\Rg{\lambda}}.
\ee
By virtue of \cite{tjac}, Lemma~2.4,
\be                  \label{der W q}
\frac{d}{d\lambda}W(\hat{\psi}_-(\lambda), \hat{\psi}_+(\lambda)) 
\Big|_{\rho_j} = - \sum_{k \in \Z} \hat{\psi}_-(\rho_j, k)
\hat{\psi}_+(\rho_j, k)
= - \frac{1}{c_j^{\pm}\gamma_{\pm, j}}.
\ee
Therefore  
\be     \label{derivative alpha}              
\frac{d}{d\lambda} \alpha(\lambda) \Big|_{\rho_j} = 
\frac {W^{\prime}(\hat{\psi}_-(\rho_j), \hat{\psi}_+(\rho_j))}{\Rg{\rho_j}} =
\frac {-1}
{c_j^{\pm}\gamma_{\pm,j} \Rg{\rho_j}}.
\ee 
 From (\ref{derivative alpha}) we obtain a connection between the left and 
right norming constants  
\be                 \label{gamma+ gamma- }
\gamma_{+,j}\gamma_{-,j} = \frac{1}
{(\alpha^{\prime}(\rho_j))^2 R_{2g+2}(\rho_j)}.
\ee
As a last preparation, we study the behavior of $\alpha(w)$ 
as $w \rightarrow 0$. 
By (\ref{B4jost}), 
\be        \label{W 0 q}
W(\psi_-(w), \psi_+(w)) = 
A_-(0) A_+(0) \tilde a w^{-1} + O(w)
\ee
and
\be
\frac {\Rg{\lambda(w)}}
{\prod_{j=1}^g(\lambda(w) - \lambda_j)} = \tilde a w^{-1} + O(1),
\ee
therefore $\alpha^{-1}(w)$ is bounded at $0$ with
\be
\alpha(0) = \prod_{j=-\infty}^{\infty}\frac{a_q(j)}{a(j)}.
\ee

We now define the {\bf scattering matrix}
\be
S(w) = \left( \begin{array}{cc} 
T(w) & R_-(w) \\
R_+(w) & T(w) \\
\end{array} \right), \qquad |w| = 1,
\ee
where $T(w) := \alpha^{-1}(w)$ and 
$R_{\pm}(w) := \alpha^{-1}(w)\beta_{\pm}(w)$ are called {\bf transmission} 
and {\bf reflection coefficients}. 
Equations (\ref{alpha beta 1}) and 
(\ref{alpha beta 2}) imply

\begin{lemma}             \label{le: S unitary}
The scattering matrix $S(w)$ is unitary.
The coefficients $T(w)$, $R_{\pm}(w)$ are bounded for $|w| = 1$, 
continuous for $|w| = 1$ except at possibly $w_l=w(E_l)$, fulfill
\bea                             \label{T R = 1}
|T(w)|^2 + |R_{\pm}(w)|^2 &=& 1, \qquad |w|=1,          \\
T(w)R_+(\overline w)+T(\overline w)R_-(w) &=& 0, \qquad |w|=1
\eea
and $\overline{T(w)} = T(\overline w)$,
$\overline{R_{\pm}(w)}=R_{\pm}(\overline w)$ for $|w|=1$

Moreover, $R^{1/2}_{2g+2}(w) T(w)^{-1}$ is continuous (in particular $T(w)$
can only vanish at $w_l$) and
\be
\label{T w_l q}
\begin{array}{l@{\quad}l}
\lim\limits_{w \rightarrow w_l} R^{1/2}_{2g+2}(w) \frac{R_{\pm}(w) + 1}{T(w)}= 0,
& w_l \neq w(\mu_j)\\
\lim\limits_{w \rightarrow w_l} R^{1/2}_{2g+2}(w) \frac{R_{\pm}(w) - 1}{T(w)}= 0,
& w_l = w(\mu_j)
\end{array}.
\ee

The transmission coefficient $T(w)$ has a meromorphic continuation to $\W$ with
simple poles at $w(\rho_j)$,
\be                                                              
\left(\res_{\rho_j}T(\lambda) \right)^2 = \gamma_{+,j}\gamma_{-,j}
R_{2g+2}(\rho_j).
\ee  
In addition, $T(z) \in \R$ as $z \in \R \backslash \sigma(H_q)$ and
\be
T(0) = \frac{1}{K_+(n,n)K_-(n,n)} = \prod_{m=-\infty}^\infty \frac{a(m)}{a_q(m)},
\ee
where $K_{\pm}(n,n)$ are the coefficients of the transformation operators.
\end{lemma}

\begin{proof}
To show (\ref{T w_l q}) we use the definition (\ref{alpha quasi}),
\[
R^{1/2}_{2g+2}(\lambda)  \frac{R_{\pm}(\lambda) + 1}{T(\lambda)} =
\prod_{j=1}^g (\lambda - \mu_j)  \big(W(\psi_-(\lambda), \psi_+(\lambda)) \mp
W(\psi_\mp(\lambda), \overline{\psi_\pm(\lambda)})\big).
\]
There are two cases to distinguish: If $\mu_j \ne E_l$ then $\psi_\pm$ are continuous
and real at $\lambda=E_l$ and the two Wronskians cancel. Otherwise,
if $\mu_j = E_l$ they are purely imaginary (by property (B2) of the Jost functions) and
the two terms are equal in the limit and add up.
\end{proof}

The sets 
\be
S_{\pm}(H) = \{R_{\pm}(w), |w| = 1; \, (\rho_j, \gamma_{\pm, j}), 
1 \leq j \leq q\}
\ee
are called left/right {\bf scattering data} for $H$.

First we want to show that the transmission coefficient can be reconstructed 
from either left or right scattering data.

Let $g(w,w_0)$ be the Green function associated with
$\W$ and let
\be
\mu(w,w_0) dw_0 =
\frac{\partial g}{\partial r} (w, r \E^{\I\theta}) \Big|_{r=1^-} 
\E^{\I\theta}d\theta,
\qquad w_0= \E^{\I\theta},
\ee
be the corresponding harmonic measure on the boundary (see, e.g., \cite{tsu}).
Since $W_0$ is simply connected, we can choose a function $h(w,v)$
such that $\hat{g}(w,w_0) = g(w,w_0) + \I h(w,w_0)$ is analytic in $W_0$.
Clearly $\hat{g}$ is only well-defined up to an imaginary constant and it will
not be analytic on $\W \backslash\{ 0\}$ in general. Similarly we can find a corresponding
$\nu(w,w_0)$ and set $\hat{\mu}(w,w_0)= \mu(w,w_0) + \I \nu(w,w_0)$. 

\begin{theorem}
Either one of the sets $S_{\pm}(H)$ determines the other and $T(w)$ via
the Poisson-Jensen type formula
\be                     \label{T(w) q}
T(w) = \exp \Bigg(
\sum_{j=1}^q  \hat{g}(w,w(\rho_j)) \Bigg)
\exp \left(\frac{1}{2}\int_{|w|=1} 
\ln(1 - |R_{\pm}(w_0)|^2) \hat{\mu}(w,w_0) dw_0 \right),
\ee
where the constant of $\hat{g}$ has to be chosen such that $T(0)>0$,
and
\[
\frac{R_-(w)}{R_+(\overline w)} = - \frac{T(w)}{T(\overline w)}, \qquad
\gamma_{+,j}\gamma_{-,j}  =  \frac{\big(\res_{\rho_j}T(\lambda)\big)^2}
{\prod_{l=0}^{2g+1}(\rho_j - E_l)}.
\] 
\end{theorem}

\begin{proof}
It suffices to prove the formula for $T(w)$, since evaluating the residua
provides $\gamma_{\pm,j}$, together with $\{\lambda_l\}$, $\{E_l\}$.
The formula for $T(w)$ holds by \cite{voza}, Theorem~1 at least when
taking absolute values. Since both sides are analytic, and have equal
absolute values, they can only differ by a constant of absolute value one.
But both sides are positive at $w=0$ and hence this constant is one.
\end{proof}

Note that neither the Blaschke factors nor the outer function in (\ref{T(w) q})
are single valued on $\W$ in general. In particular, the eigenvalues cannot
be chosen arbitrarily, which was first observed in \cite{kumi}.

\section{The Gel'fand-Levitan-Marchenko equations}
\label{secGLM}

In this section we want to derive a procedure which allows the
reconstruction of the Jacobi operator $H$ with asymptotically quasi-periodic 
coefficients from its scattering data $S_{\pm}(H)$. This will be achieved
by deriving an equation for $K_{\pm}(n, m)$ which is generally known as 
Gel'fand-Levitan-Marchenko equation.

Since $K_{\pm}(n, m)$ are essentially the Fourier coefficients of the
Jost solutions $\psi_{\pm}(w, n)$ we compute the Fourier
coefficients of the scattering relations (\ref{scat rel q}). Therefore 
we multiply 
\be                                              \label{T psi}
T(w) \psi_{\mp}(w, n) = R_{\pm}(w) \psi_{\pm}(w, n) +
\overline{\psi_{\pm}(w, n)} 
\ee
by $(2 \pi \I)^{-1} \psi_{q, \pm}(w, m) d\omega$, where 
$\pm m \geq \pm n$, and integrate around the unit circle. 
First we evaluate the right hand side of (\ref{T psi}) using (\ref{de:K qp})
\bea                                            \label{GLM 1 q}
\frac{1}{2\pi \I} \int_{|w|=1} \overline{\psi_+(w,n)} 
\psi_{q, +}(w, m) d\omega(w) &=& K_+(n,m),       \\ \nonumber
\frac{1}{2\pi \I} \int_{|w|=1} R_+(w) \psi_+(w, n)
\psi_{q, +}(w, m) d\omega(w) &=& \sum_{l=n}^{\infty}K_+(n,l)
\tilde F^+(l, m),
\eea
where
\be       \label{fourier R q}
\tilde F^+(l, m) = \frac{1}{2\pi \I}\int_{|w| = 1}
R_+(w) \psi_{q, +}(w, l) \psi_{q, +}(w, m) d\omega(w).
\ee
Note that $\tilde F^+(l, m) = \tilde F^+(m, l)$ is real.

To evaluate the left hand side of (\ref{T psi}) we use the residue theorem.
The only poles are at the eigenvalues and at $0$ if $n=m$, hence
\bea     \nonumber
&&\frac{1}{2\pi \I} \int_{|w|=1}
T(w) \psi_-(w, n)\psi_{q, +}(w, m)d\omega(w)  \\ \nonumber
&& = \frac{\delta(n,m)}{K_+(n,n)} 
+ \sum_{j=1}^q  \res_{\rho_j} 
\Big(\frac{T(\lambda) \hat{\psi}_-(\lambda, n)
\hat{\psi}_{q, +}(\lambda, m)}
{\Rg{\lambda}}\Big).
\eea
Here $\delta(n,m)$ is one for $m=n$ and zero else. By (\ref{derivative alpha})
the residua at the eigenvalues are given by
\be
\res_{\rho_j}
\Big(\frac{T(\lambda) \hat{\psi}_-(\lambda, n)
\hat{\psi}_{q, +}(\lambda, m)}
{\Rg{\lambda}}\Big) = - \gamma_{+,j}  \hat{\psi}_+(\rho_j, n)
\hat{\psi}_{q, +}(\rho_j, m).
\ee
Collecting all terms yields
\be
K_{\pm}(n, m) + \sum_{l=n}^{\pm \infty}K_{\pm}(n,l)\tilde F^{\pm}(l, m)
= \frac{\delta(n,m)}{K_{\pm}(n,n)} - \sum_{j=1}^q
\gamma_{\pm,j} \hat{\psi}_{\pm}(\rho_j, n) \hat{\psi}_{q, \pm}(\rho_j, m)
\ee 
and we have thus proved the following result.

\begin{theorem}
The kernel $K_{\pm}(n,m)$ of the transformation operator satisfies
the {\bf Gel'fand-Levitan-Marchenko equation}
\be   \label{glm1 q}
K_{\pm}(n,m) + \sum_{l=n}^{\pm \infty}K_{\pm}(n,l)F^{\pm}(l,m) = 
\frac{\delta(n,m)}{K_{\pm}(n,n)},  \qquad \pm m \geq \pm n,
\ee
where
\be          \label{glm2 q}
F^{\pm}(l,m) =  \tilde F^{\pm}(l, m) + \sum_{j=1}^q \gamma_{\pm,j} 
\hat{\psi}_{q, \pm}(\rho_j, l) \hat{\psi}_{q, \pm}(\rho_j, m).
\ee
\end{theorem}

Defining the {\bf Gel'fand-Levitan-Marchenko operator}   
\be
\mathcal{F}^{\pm}_n f(j) = \sum_{l=0}^\infty F^{\pm}(n \pm l, n \pm j)f(l),
\qquad f \in \ell^2(\N_0, \C), 
\ee
yields that the Gel'fand-Levitan-Marchenko equation is equal to  
\be       \label{marchenko 2 q}              
(1 + \mathcal{F}^{\pm}_n) K_{\pm}(n, n \pm .) = (K_{\pm}(n, n))^{-1} \delta_0.
\ee

Our next aim is to study the Gel'fand-Levitan-Marchenko operator 
$\mathcal{F}^{\pm}_n$ in more detail. The structure of the 
Gel'fand-Levitan-Marchenko equation suggests that the estimate 
(\ref{estimate K q}) for $K_{\pm}(n, m)$ should imply a similar estimate 
for $F^{\pm}(n,m)$.   

\begin{lemma}          \label{le: estimate F q}
\be             \label{estimate F q}
|F^{\pm}(n, m)| \leq C \sum_{j=[\frac{n+m}{2}] \pm 1}^{\pm \infty} 
\Big(|a(j) - a_q(j)| + |b(j) - b_q(j)| \Big), 
\ee
where the constant $C$ is of the same nature as in (\ref{estimate K q}).
\end{lemma}
 
\begin{proof}
We abbreviate the estimate (\ref{estimate K q}) for $K_+(n,m)$ by 
\be       \label{est K_+ q}
|K_+(n,m)| \leq C\, C_+(n+m),
\ee
where
\[
C_+(n+m) = \sum_{j=[\frac{n+ m}{2}]+1}^{\infty} c(j), \qquad
c(j) = |a(j) - a_q(j)| + |b(j) - b_q(j)|.
\]
Note that $C_+(n+1) \leq C_+(n)$. Moreover, $C_+(n) \in \ell_+^1(\Z)$ since 
the summation by parts formula (e.g.\
\cite{tjac}, (1.18))
\be
\sum_{m=n}^N g(m)(f(m+1)-f(m)) = g(N)f(N+1)-g(n-1)f(n) + \sum_{m=n}^N
(g(m-1)-g(m))f(m)
\ee
implies for $g(m)=m$, $f(m)=C_+(m)$ that
\be \label{nc(n)}
\sum_{m=n}^{\infty}m\, c(m) = (n-1)C_+(n) + \sum_{m=n}^{\infty}C_+(m),
\ee
where we used   $\lim_{n \rightarrow \infty} n\, C_+(n+1) 
\leq \lim_{n \rightarrow \infty} \sum_{m=n}^{\infty} m\, c(m)=0$.
Solving the GLM-equation (\ref{glm1 q}) for $F^+(n, m)$, $m > n$, we obtain 
\bea              \nonumber
|F^+(n, m)| &\leq& \frac{1}{K_+(n,n)}
\left(|K_+(n,m)| + \sum_{l=n+1}^{\infty}
\left|K_+(n, l) F^+(l, m)\right|\right)    \\
&\leq& C_1(n) \left(C_+(n+m) + 
\sum_{l=n+1}^{\infty}C_+(n+l)           \nonumber
\left|F^+(l, m)\right|\right),         
\eea
where $C_1(n) = C\, |K_+(n,n)|^{-1} \rightarrow C$ for 
$n \rightarrow \infty$ (see (\ref{K(n,n)})). 
For $n$ large enough, i.e.\ $C_1(n) C_+(2n) < 1$, we apply the discrete 
Gronwall-type inequality \cite{tjac}, Lemma~10.8,
\bea                          \nonumber
|F^+(n, m)| &\leq&  C_1(n) \left(C_+(n + m) + \sum_{l=n+1}^{\infty}
\frac{C_1(l) C_+(l+m)C_+(n+l)}
{\prod_{k=n+1}^l (1- C_1(k) C_+(n+k))}\right) \\ \label{est F^+ q}
&\leq& C_1(n) C_+(n + m) \left(1 + \sum_{l=n+1}^{ \infty}         
\frac{C_1(k) C_+(n+l)}{\prod_{k=n+1}^l (1-C_1(n) C_+(n+k))}\right),
\eea 
which finishes the proof.
\end{proof}

Furthermore,

\begin{lemma}       \label{le: estimate F dif}
Let $F^{\pm}(n,m)$ be solutions of the Gel'fand-Levitan-Marchenko equation.
Then
\bea           \label{sum F+- q}                          
\sum_{n=n_0}^{\pm \infty}|n| \left|F^{\pm}(n,n) 
- F^{\pm}(n \pm1, n \pm 1)\right|
&<& \infty,              \\      \label{sum F+ q}                   
\sum_{n=n_0}^{\pm \infty}|n| \left|a_q(n)F^{\pm}(n,n+1) 
- a_q(n-1)F^{\pm}(n -1, n)\right|
&<& \infty. 
\eea
\end{lemma}

\begin{proof}
We first prove (\ref{sum F+ q}) for $F^+$. Lemma~\ref{coro K+K- q} implies  
\be          \label{b-b_q}                
b(n) - b_q(n) = a_q(n)\kappa_{+,1}(n)-a_q(n-1)\kappa_{+,1}(n-1),   
\ee
where 
\be
\kappa_{+,j}(n) := \kappa_{+}(n, n+j) := \frac{K_+(n, n+j)} {K_+(n,n)}.
\ee 
Abbreviate $F_j^+(n) := F^+(n+j, n)$. With this notation, the  
GLM-equation (\ref{glm1 q}) reads
\be
\kappa_{+,l}(n) + F_l^+(n) + \sum_{j=1}^{\infty} \kappa_{+,j}(n) F_{j-l}^+(n+l) = 
\frac{\delta(l,0)}{K_+(n,n)^2},
\qquad l \geq 0.
\ee
Insert the GLM-equation for $F^+(n, n+1)$, $F^+(n-1, n)$ (recall 
$F^+(n,m)=F^+(m,n)$)
\bea   \nonumber
\lefteqn{a_q(n)F_1^+(n) - a_q(n-1)F_1^+(n-1)}    \\  \nonumber 
&=& - a_q(n)\kappa_{+,1}(n) + a_q(n-1)\kappa_{+,1}(n-1)
\\
&& - \sum_{j=1}^{\infty}\Big( a_q(n) \kappa_{+,j}(n)
F_{j-1}^+(n+1) - a_q(n-1) \kappa_{+,j}(n-1) F_{j-1}^+(n)\Big).
\eea
Since $- a_q(n)\kappa_{+,1}(n) + a_q(n-1)\kappa_{+,1}(n-1)= b_q(n) - b(n)$
the only interesting part is the sum. For $N$, $J < \infty$,
\bea     \label{sum F+ 1}                                               \nonumber
\lefteqn{\sum_{n=n_0}^N n \sum_{j=1}^J \Big( a_q(n) \kappa_{+,j}(n)
F_{j-1}^+(n+1) - a_q(n-1) \kappa_{+,j}(n-1) F_{j-1}^+(n)\Big)  }
\\ \nonumber
&=& \sum_{j=1}^J \sum_{n=n_0}^N n \Big( a_q(n) \kappa_{+,j}(n)
F_{j-1}^+(n+1) - a_q(n-1) \kappa_{+,j}(n-1) F_{j-1}^+(n) \Big)
\\ \nonumber
&=& \sum_{j=1}^J \Big( N a_q(N) \kappa_{+,j}(N) F_{j-1}^+(N+1) -
(n_0-1) a_q(n_0-1) \kappa_{+,j}(n_0-1) F_{j-1}^+(n_0)
\\  
&& \qquad + \sum_{n=n_0}^N (-1) a_q(n-1) \kappa_{+,j}(n-1) F_{j-1}^+(n) \Big),
\eea
where we used the summation by parts.
Estimates (\ref{est K_+ q}), (\ref{est F^+ q}) imply for the first summand
\bea \nonumber
\Big| \sum_{j=1}^J N a_q(N) \kappa_{+,j}(N) F_{j-1}^+(N+1) \Big| &\leq&
\sum_{j=1}^J |N| a_q(N) \tilde C C_+(2N + j) C_+(2N +j+1)
\\ \nonumber
&\leq& |N| a_q(N) \hat C C_+(2N + 1),
\eea
which holds uniformly in $J$, and (compare (\ref{nc(n)}))
\be
\lim_{N \rightarrow \infty} N a_q(N) \hat C C_+(2N + 1) = 0.
\ee
Moreover,
\bea \nonumber
\lefteqn{ \lim_{N, J \rightarrow \infty}
\Big| \sum_{j=1}^J \sum_{n=n_0}^N a_q(n-1) \kappa_{+,j}(n-1) F_{j-1}^+(n) \Big|}
\\ \nonumber &\leq& 
\lim_{N, J \rightarrow \infty}
\sum_{j=1}^J \sum_{n=n_0}^N \Big| a_q(n-1) \kappa_{+,j}(n-1) F_{j-1}^+(n) \Big|
\\ \nonumber
&\leq& \sum_{j=1}^{\infty} \sum_{n=n_0}^{\infty} a_q(n-1) \tilde C C_+(2n + j) C_+(2n +j+1)
< \infty.
\eea
Therefore 
$|n||a_q(n)F^+(n,n+1) - a_q(n-1)F^+(n -1, n)| \in \ell^1_+(\Z)$
as desired.
To apply Lemma \ref{coro K+K- q} for $F^-$ use the 
symmetry property $F^-(n,m)=F^-(m,n)$. 
For (\ref{sum F+- q}), inserting the GLM-equation yields  
\bea  \nonumber                    
\lefteqn{F^+(n,n) - F^+(n+1, n + 1) = K_+^{-2}(n, n) - K_+^{-2}(n+1, n+1)} 
\\   \nonumber
&& +  \sum_{j=1}^{\infty} \Big( \kappa_{+,j}(n+1) F_j^+(n+1) -
\kappa_{+,j}(n) F_j^+(n)\Big).   
\eea
By (\ref{K(n,n)}),  
\bea                    \nonumber
\left| K_+^{-2}(n, n) - K_+^{-2}(n+1, n+1)\right| &\leq&
\frac{|a(n)+a_q(n)|}{a(n)^2}\prod_{j=n+1}^{\infty} \frac{a(j)^2}{a_q(j)^2} 
|a(n) - a_q(n)| \\  &\leq& C |a(n)-a_q(n)|,               \label{F diff 4}
\eea
and the same considerations as above imply (\ref{sum F+- q}).
\end{proof}

\begin{remark}
The Gel'fand-Levitan-Marchenko equation
is symmetric in $K_{\pm}(n,m)$ and $F^{\pm}(n,m)$, therefore
we can invert the analysis done in Lemma \ref{le: estimate F dif} and obtain
estimates for $K_{\pm}(n,m)$ starting with an analogue of estimate 
(\ref{estimate F q}) for $F^{\pm}(n,m)$ and the estimates 
(\ref{sum F+- q}), (\ref{sum F+ q}) (c.f.\ Lemma \ref{invs1 q}).
\end{remark}

\begin{theorem} \label{thmglm}
For $n \in \Z$, the Gel'fand-Levitan-Marchenko operator
$\mathcal{F}^{\pm}_n: \ell^2 \to \ell^2$ is Hilbert-Schmidt.
Moreover, $1+\mathcal{F}^{\pm}_n$ is positive and hence invertible.

In particular, the Gel'fand-Levitan-Marchenko equation (\ref{marchenko 2 q}) 
has a unique solution and $S_+(H)$ or $S_-(H)$ uniquely determine $H$.
\end{theorem}

\begin{proof}
That $\mathcal{F}^{\pm}_n$ is Hilbert-Schmidt is a straight-forward consequence of our estimate Lemma~\ref{le: estimate F q}.

Let $f \in \ell^2(\N_0)$ be real
(which is no restriction since $F^+(n,l)$ is real and the 
real and imaginary part of (\ref{unique 1}) could be treated separately) 
and abbreviate 
$f_n(w) = \sum_{j=0}^{\infty} f(j)   \psi_{q,+}(w, n+j)$.
Then
\bea      \nonumber
&&\sum_{j=0}^{\infty} f(j) \mathcal{F}_n^+ f(j) =
\sum_{j=0}^{\infty} f(j) \sum_{l=0}^{\infty} 
F^+(n+j, n+l) f(l)      \\ \nonumber
&& = \frac{1}{2 \pi i} \int_{|w|=1} R_+(w) \sum_{j,l =0}^{\infty}
f(j)   \psi_{q,+}(w, n+j)   \psi_{q,+}(w, n+l)      f(l)
\, d\omega(w)     \\ \nonumber
&& \qquad + \sum_{k=1}^q \sum_{j,l =0}^{\infty} f(j) \gamma_{+,k}
\hat \psi_{q,+}(\rho_k, n+j) \hat \psi_{q,+}(\rho_k, n+l) f(l)
\\ \nonumber
&&= \frac{1}{2 \pi i} \int_{|w|=1} R_+(w) f_n(\overline w)
f_n(w)\, d\omega(w) + \sum_{k=1}^q \gamma_{+,k} 
|\hat{f}_n(\rho_k)|^2  \\ \label{unique 1}
&&= \frac{1}{2 \pi i} \int_{|w|=1} \tilde R_+(w) 
|f_n(w)|^2 d\omega(w)
+ \sum_{k=1}^q \gamma_{+,k} |\hat{f}_n(\rho_k)|^2,
\eea
where $\tilde R_+(w) = R_+(w) f_n(w)\, \big(\overline{f_n(w)}\big)^{-1}$
with $|\tilde R_+(w)| = |R_+(w)|$ and $\hat{f}_n(w) =
\sum_{j=0}^{\infty} f(j)  \hat{\psi}_{q,+}(w, n+j)$.  The integral over the imaginary
part vanishes since $\overline{\tilde R_+(w)} = \tilde R_+(\overline w)$
and we replace the real part by
\be     \nonumber
\mbox{Re}(\tilde R_+(w)) = \frac{1}{2}\big(|1+\tilde R_+(w)|^2 - 1 - 
|\tilde R_+(w)|^2 \big) = \frac{1}{2}\big(|1+\tilde R_+(w)|^2 + 
|T(w)|^2 \big) - 1,
\ee
(recall $|\tilde R_+(w)|^2 + |T(w)|^2 = 1$). 
This yields using $\sum |f(j)|^2 = \frac{1}{2 \pi i}\int_{|w|=1} 
|f_n(w)|^2 d\omega$ 
\bea     \nonumber
\sum_{j=0}^{\infty} \overline{f(j)}(1 + \mathcal{F}_n^+) f(j) &=&
\sum_{k=1}^q \gamma_{+,k} |\hat{f}_n(\rho_k)|^2  \\ 
&& + \frac{1}{4 \pi i}\int_{|w|=1} 
\big(|1+\tilde R_+(w)|^2 + |T(w)|^2 \big) |f_n(w)|^2 d\omega(w)
\eea
which establishes $1 + \mathcal{F}_n^+ \geq 0$. According to 
Lemma~\ref{le: S unitary}, $|T(w)|^2 > 0$ a.e., therefore
$-1$ is not an eigenvalue and $1 + \mathcal{F}_n^+ \geq \epsilon_n$
for some $\epsilon_n > 0$.
\end{proof}

To finish the direct scattering step for the Jacobi operator $H$ with
asymptotically quasi-periodic coefficients we summarize the properties of
the scattering data $S_{\pm}(H)$.

\begin{hypo}             \label{hypo scat q}
The scattering data
\be
S_{\pm}(H) = \{R_{\pm}(w), |w|=1; (\rho_j, \gamma_{\pm, j}), 1 \leq j
\leq q\}
\ee
satisfy the following conditions.

(i). The reflection coefficients $R_{\pm}(w)$ are continuous except
possibly at $w_l=w(E_l)$ and fulfill
\be
\overline{R_{\pm}(w)} = R_{\pm}(\overline w).
\ee
Moreover, $|R_\pm(w)|<1$ for $w\neq w_l$ and
\be \label{hypo T}
1-|R_\pm(w)|^2 \ge C \prod_{l=0}^{2g+1} |w- w_l|^2.
\ee
The Fourier coefficients
\be
\tilde F^\pm(l, m) = \frac{1}{2\pi \I}\int_{|w| = 1}
R_\pm(w) \psi_{q, \pm}(w, l) \psi_{q, \pm}(w, m) d\omega(w)
\ee
satisfy  
\bea         \nonumber                      
&&|\tilde F^{\pm}(n, m)| \leq \sum_{j=n+m}^{\pm \infty} q(j), 
\qquad q(j) \geq 0, \qquad  |j| q(j) \in \ell^1(\Z), \\         \nonumber
&&\sum_{n = n_0}^{\pm \infty}|n| \big|\tilde F^{\pm}(n,n) - 
\tilde F^{\pm}(n \pm 1, n \pm 1)\big| < \infty,    \\ \nonumber
&&\sum_{n = n_0}^{\pm \infty}|n| \big|a_q(n) \tilde F^{\pm}(n,n+1) - 
a_q(n-1) \tilde F^{\pm}(n - 1, n)\big| < \infty.
\eea

(ii). The values $\rho_j  \in \R\backslash \sigma(H_q)$, 
$1 \leq j \leq q$, are distinct and the 
norming constants $\gamma_{\pm, j}$, $1 \leq j \leq q$, 
are positive.

(iii). $T(w)$ defined via equation (\ref{T(w) q})
extends to a single valued function on $\W$ (i.e., it has equal values on the
corresponding slits).

(iv). 
Transmission and reflection coefficients satisfy
and satisfies
\be  \label{lim T}
\begin{array}{l@{\qquad}l}
\lim\limits_{w \rightarrow w_l} (w-w_l) \frac{R_{\pm}(w) + 1}{T(w)} = 0,
&  w_l \neq w(\mu_j) \\
\lim\limits_{w \rightarrow w_l} (w-w_l) \frac{R_{\pm}(w) - 1}{T(w)} = 0
&   w_l = w(\mu_j)
\end{array}.
\ee
and the consistency conditions
\[              
\frac{R_-(w)}{R_+(\overline w)} = - \frac{T(w)}{T(\overline w)}, \qquad
\gamma_{+, j}\, \gamma_{-, j} =  
\frac{\big(\res_{\rho_j}T(\lambda)\big)^2}
{\prod_{l=0}^{2g+1}(\rho_j - E_l)}.
\]
\end{hypo} 

\begin{remark}
Note that (\ref{hypo T}) implies that $\ln(1-|R_\pm(w)|^2)$ is integrable
and ensures that (\ref{T(w) q}) is well-defined, at least as a multi valued
function. Condition (iii), which is void in the constant background case,
shows that the the reflection coefficient and eigenvalues cannot be
chosen independent of each other.
\end{remark}

\section{Inverse scattering theory}
\label{secINV}

In this section we want to invert the process of scattering theory, that
is, we want to reconstruct the operator $H$ from a given set $S_\pm$ and 
a given quasi-periodic Jacobi operator $H_q$.

If $S_{\pm}$ (satisfying H.\ref{hypo scat q} (i)--(ii)) and $H_q$ are known, we
can construct $F^{\pm}(l,m)$ via formula (\ref{glm2 q}) and thus derive the
Gel'fand-Levitan-Marchenko  equation, which has a unique solution by
Theorem~\ref{thmglm}. This solution 
\bea \nn
K_\pm(n,n) &=& \spr{\delta_0}{(1 + \mathcal{F}^{\pm}_n)^{-1} \delta_0}^{1/2}\\
K_\pm(n,n\pm j) &=& \frac{1}{K_\pm(n,n)} \spr{\delta_j}{(1 + \mathcal{F}^{\pm}_n)^{-1} \delta_0}
\eea
is the kernel of the transformation operator. Since $1 + \mathcal{F}^{\pm}_n$ is
positive, $K_\pm(n,n)$ is positive and we can set in accordance with 
Lemma~\ref{coro K+K- q} 
\bea                             \label{a+a- q}
a_+(n) &=& a_q(n)\frac{K_+(n+1, n+1)}{K_+(n,n)}, \\             \nonumber
a_-(n) &=& a_q(n)\frac{K_-(n, n)}{K_-(n+1,n+1)}, \\             \nonumber
b_+(n) &=& b_q(n) + a_q(n)\frac{K_+(n, n+1)}{K_+(n,n)} - a_q(n-1)
\frac{K_+(n-1, n)}{K_+(n-1,n-1)}, \\  \nonumber           
b_-(n) &=& b_q(n) + a_q(n-1)\frac{K_-(n, n-1)}{K_-(n,n)} - a_q(n)
\frac{K_-(n+1, n)}{K_-(n+1,n+1)}.
\eea
Let $H_+$, $H_-$ be the associated Jacobi operators.

\begin{lemma}             \label{invs1 q}
Suppose a given set $S_{\pm}$ satisfies H.\ref{hypo scat q} (i)--(ii). Then the
sequences defined in (\ref{a+a- q}) satisfy $n |a_{\pm}(n)-a_q(n)|$,
$n|b_{\pm}(n)-b_q(n)| \in \ell^1_\pm(\N)$. 

Moreover, $\psi_{\pm}(\lambda, n) = \sum_{m=n}^{\pm \infty}
K_{\pm}(n,m)\psi_{q,\pm}(\lambda, m)$, where $K_{\pm}(n,m)$ is the 
solution of the Gel'fand-Levitan-Marchenko equation, satisfies 
$\tau_{\pm} \psi_{\pm} = \lambda \psi_{\pm}$.
\end{lemma}

\begin{proof}
We only prove the statements for the "+" case.
Define $F^+(n,m)$ by (c.f.\ (\ref{glm2 q}))
\[          
F^+(l,m) =  \tilde F^+(l, m) + \sum_{j=1}^q \gamma_{+,j} 
\hat \psi_{q, +}(\rho_j, l) \hat \psi_{q, +}(\rho_j, m).
\]
Hypothesis H.\ref{hypo scat q} (i) implies
\bea             \label{inv 1}
&&|F^+(n, m)| \leq C \sum_{j=n+m}^{\infty} q(j) =: C_+(n+m), \\
&&\sum_{n = n_0}^{\infty}|n| \big| F^+(n,n) -  \label{inv F 2}
F^+(n + 1, n + 1)\big| < \infty,    \\ \label{inv F 1}
&&\sum_{n = n_0}^{\infty}|n| \big|a_q(n) F^+(n,n+1) - 
a_q(n-1) F^+(n - 1, n)\big| < \infty,
\eea
since $\hat \psi_{q, +}(\rho_j, n)$ decay exponentially as $n \rightarrow \infty$
and $\sum_{j} \gamma_{+,j} \hat \psi_{q, +}(\rho_j, .) \hat \psi_{q, +}(\rho_j, .)$
form a telescopic sum. Note that $C_+(n+1) < C_+(n)$. 

Set $\kappa_+(n,m):= K_+(n,m)K_+(n,n)^{-1}$. Then as in the proof of
Lemma~\ref{le: estimate F q} we obtain
\be     \label{inv 2}                
|\kappa_+(n,m)| \leq  C_+(n+m)(1 + O(1)).
\ee
Now we have all estimates at our disposal to prove 
$n|b_+(n)-b_q(n)| \in \ell^1(\N)$. By definition (c.f.\ (\ref{a+a- q})),
\be 
b_+(n) - b_q(n) = a_q(n)\kappa_+(n, n+1) - a_q(n-1)\kappa_+(n-1, n).
\ee
We insert the GLM-equation for $\kappa_+(n, n+1)$, $\kappa_+(n-1, n)$ and use
estimate (\ref{inv F 1}), the summation by parts formula, and estimates (\ref{inv 1}), 
(\ref{inv 2}) in the same way as in Lemma~\ref{le: estimate F dif}. 
Similarly using (\ref{inv F 2}) we see 
\be
\sum_{n=n_0}^{\infty}|n|\left|\frac{1}{K_+^2(n,n)} 
- \frac{1}{K_+^2(n+1,n+1)}\right| < \infty.
\ee
Equation (\ref{a+a- q}) yields
\[
\left|\frac{1}{K_+^2(n,n)} - \frac{1}{K_+^2(n+1,n+1)}\right| =
\frac{1}{a_q(n)^2}\Big(\prod_{j=n+1}^{\infty}\frac{a_+(j)^2}
{a_q(j)^2}\Big)|a_+(n)^2 - a_q(n)^2|.
\]
The product converges and therefore $|n||a_+(n)^2 - a_q(n)^2| \in \ell^1(\N)$.

Next we consider $\psi_+(\lambda, n)$. Abbreviate
\bea         \label{Delta K q}
&& (\Delta K_+)(n,m) = a_q(n-1)\kappa_+(n-1,m) +
a_+^2(n)a_q^{-1}(n)\kappa_+(n+1,m)  \\           \nonumber
&& - a_q(m-1)\kappa_+(n,m-1) -
a_q(m)\kappa_+(n,m+1) + (b_+(n)-b_q(m))\kappa_+(n,m).  
\eea
$\Delta K_+ = 0$ is equivalent to the operator equality 
$H_+K_+=K_+H_q$, which in turn implies that 
$\psi_+(\lambda,n)$ satisfies 
$H_+ \psi_+ = \lambda \psi_+$ 
\be
H_+ \psi_+ = H_+K_+ \psi_{q,+} =
K_+H_q \psi_{q,+} = K_+\lambda \psi_{q,+} =
\lambda K_+ \psi_{q,+} = \lambda \psi_+.
\ee
To show that $\Delta K_+ = 0$ we insert
the GLM-equation into (\ref{Delta K q}) and obtain 
\be \label{delta K}
(\Delta K_+)(n,m) + \sum_{l=n+1}^{  \infty}
(\Delta K_+)(n,l)F^+(l,m) = 0, \quad m > n+1.
\ee
In the calculations we used
\bea            \nonumber
&&a_q(n-1)F^+(n-1,m) + b_q(n)F^+(n,m) + a_q(n)F^+(n+1,m) = \\
&&a_q(m-1)F^+(n,m-1) + b_q(m)F^+(n,m)+ a_q(m)F^+(n,m+1)
\nonumber
\eea
which follows from (\ref{glm2 q}).   
By Theorem~\ref{thmglm} equation (\ref{delta K}) has only the 
trivial solution $\Delta K_+=0$ and hence the proof is complete.
\end{proof}

Now we can prove the main result of this section.

\begin{theorem}
Hypothesis H.\ref{hypo scat q} is necessary and sufficient for a sets $S_\pm$
to be the left/right scattering data of a unique Jacobi operator $H$ associated with
sequences $a$, $b$ satisfying H.\ref{hypo:quasi}.
\end{theorem}

\begin{proof}
Necessity has been established in the previous section. By Lemma~\ref{invs1 q}, we
know existence of sequences $a_\pm$, $b_\pm$ and corresponding solutions
$\psi_\pm(w,n)$ associated with $S_+$ (or $S_-$). Hence it remains to establish
$a_+(n)=a_-(n)$ and $b_+(n)=b_-(n)$.

Consider the following part of the GLM-equation  
\be
\Phi_+(n, .) := \sum_{l=n}^{\infty}K_+(n,l)
{\tilde F}^+(l,.) \in \ell_+^1(\Z).
\ee
Then by use of (\ref{GLM 1 q}) and Lemma \ref{lem orth},
\bea             \nonumber
\lefteqn{\sum_{m \in \Z}\Phi_+(n, m) \psi_{q,-}(w, m) =
\sum_{m \in \Z}\bigg(\sum_{l=n}^{\infty}K_+(n,l)
{\tilde F}^+(l,m)\bigg) \psi_{q,-}(w, m)}  \\         \nonumber
&=& \sum_{m \in \Z} \bigg( \frac{1}{2\pi \I} \int_{|w|=1} 
R_+(w) \psi_+(w, n) \psi_{q, +}(w, m)
d\omega(w)\bigg) \psi_{q,-}(w, m) \\             \nonumber
&=& \sum_{m \in \Z} 
\langle \psi_{q,-}(w, m), R_+(w) \psi_+(w, n)
\rangle \psi_{q,-}(w, m) \\ 
&=& R_+(w) \psi_+(w, n).
\eea
On the other hand, inserting the GLM-equation 
yields for $|w|=1$
\bea
\lefteqn{\sum_{m \in \Z}\Phi_+(n, m) \psi_{q, -}(w, m) =} \nonumber \\
&=& \sum_{m = - \infty}^{n-1}\Phi_+(n, m) \psi_{q,-}(w, m) 
+ \sum_{m = n}^{\infty}  \Big[       \nonumber
\delta(n, m) K_+^{-1}(n, n) -  K_+(n, m) 
\\      \nonumber
&& - \sum_{l=n}^{\infty} K_+(n, l)   
\sum_{j=1}^q \gamma_{+,j} \hat \psi_{q, +}(\rho_j,l) 
\hat \psi_{q, +}(\rho_j,m)\Big] \psi_{q,-}(w, m) \\
&=&                       \nonumber
\sum_{m = - \infty}^{n-1}\Phi_+(n, m) \psi_{q,-}(w, m)  
+ \psi_{q,-}(w,n)K_+^{-1}(n, n) 
- \overline{\psi_+(w, n)}\\
&& - \sum_{j=1}^q \gamma_{+, j}       
\hat \psi_+(\rho_j, n) 
\sum_{m = n}^{\infty}\hat \psi_{q,+}(\rho_j, m) \psi_{q,-}(w, m),
\eea
(recall the definition of $\hat \psi_{q,\pm}$ from (\ref{def eigen}))
and therefore
\be                  \label{h3 q}
T(w)h_-(w, n) = \overline{\psi_+(w, n)} 
+ R_+(w) \psi_+(w, n), \qquad |w|=1,
\ee
where
\bea \nonumber         \label{h_mp(w, n)}
h_-(w, n) &=& \frac{\psi_{q,-}(w, n)}{T(w)}\bigg( 
\frac{1}{K_+(n, n)} + 
\sum_{m = - \infty}^{n-1}\Phi_+(n, m) 
\frac{\psi_{q,-}(w, m)}{\psi_{q,-}(w, n)}  \\
&&+ \sum_{j=1}^q \gamma_{+, j}       
\hat \psi_+(\rho_j, n)  
\frac{W_{n-1}(\hat \psi_{q,+}(\rho_j), \psi_{q,-}(w))}
{\psi_{q,-}(w, n) (\lam(w) - \rho_j)}\bigg),             
\eea
since Green's formula (\cite{tjac}, eq. (1.20)) implies for $\lambda \in 
\sigma(H_q)$  
\be     \label{green h+-}                                \nonumber
(\lambda - \rho_j)\sum_{m = n}^{\infty}
\hat \psi_{q,+}(\rho_j, m) \psi_{q,-}(\lambda, m)
= - W_{n-1}(\hat \psi_{q,+}(\rho_j), \psi_{q,-}(\lambda)).
\ee     
Similarly, we obtain
\bea \nonumber                           \label{h_p(w, n)}
h_+(w, n) &=& \frac{\psi_{q,+}(w, n)}{T(w)}\bigg( 
\frac{1}{K_-(n, n)} + 
\sum_{m = n+1}^{\infty} \Phi_-(n, m) 
\frac{\psi_{q,+}(w, m)}{\psi_{q,+}(w, n)}  \\
&&- \sum_{j=1}^q \gamma_{-, j}       
\hat \psi_-(\rho_j, n) 
\frac{W_n(\hat \psi_{q,-}(\rho_j), \psi_{q,+}(w))}
{\psi_{q,+}(w, n)(\lam(w) - \rho_j)}\bigg)             
\eea
with
\[
\Phi_-(n, m) = \sum_{l= - \infty}^n K_-(n,l)
{\tilde F}^-(l, m).
\]
For $n \in \Z$, $|w|=1$, we see that 
$h_{\mp}(w^{-1}, n)=\overline{h_{\mp}(w, n) }$, since $K_{\pm}(n, m)$ and 
$\Phi_{\pm}(n, m)$ are real. The functions $h_{\mp}(w, n)$ are continuous
for $|w|=1$, $w \neq w(E_j)$, since $T^{-1}(w)$ is
continuous on this set by the Poisson-Jensen formula (\ref{T(w) q})
($|R_{\pm}(w)|<1$ for $w \neq w(E_j)$ by H.\ref{hypo scat q} (i)) 
and $\psi_{q,\mp}(w, m)$ are continuous on 
$\partial W\backslash\{w(\mu_k)\}$. The functions $h_{\mp}(w, n)$
have a meromorphic continuation to $\W \backslash \{0\}$
with the only possible poles at $w(\rho_j)$ and $w(\mu_j)$.
At $w(\rho_j)$ there are no poles, due to the zeros of $T^{-1}(w)$ 
at $w(\rho_j)$. For $w=w(\mu_j)$ we have the same type of singularity as
$\psi_{q,\pm}$. In summary, $h_{\pm}(w, n)$ have simple poles at
$w(\mu_j)$ and are continuous at the boundary except possibly at $w(E_j)$.

To study the behavior of $h_{\pm}(w, n)$ as $w \rightarrow 0$, 
we recall $z^{-1}= - w/\ti{a}\, (1+ O(w))$. Then
\bea                     \nonumber
\lefteqn{\big(\frac{w}{\ti{a}} + O(w^2)\big)
W_{n-1}(\hat \psi_{q,+}(\rho_j), \psi_{q,-}(w))} \\ 	\nonumber
&=& \frac{(-1)^n \ti{a}^{n-1}}{\prod_{j=0}^{n-2}a_q(j)} w^{-n+1}
(\hat \psi_{q,+}(\rho_j, n-1) + O(w)),  \\ 	\nonumber
\lefteqn{\big(\frac{-w}{\ti{a}} + O(w^2)\big)
W_n(\hat \psi_{q,-}(\rho_j), \psi_{q,+}(w))} \\ \nonumber 
&=& \frac{(-1)^n \prod_{j=0}^n a_q(j)}{\tilde a^{n+1}} w^{n+1}
(\hat \psi_{q,-}(\rho_j, n+1) + O(w)),
\eea
and property (B4) implies  
\be
\sum_{m = n \mp 1}^{\mp \infty}\Phi_{\pm}(n, m) \psi_{q,\mp}(w, m)
\psi_{q,\mp}^{-1}(w, n) = O(w), \qquad w \rightarrow 0.
\ee
We conclude that
\be
\lim_{w\rightarrow 0} h_{\mp}(w,n) \psi_{q, \pm}(w,n) = \frac{1}
{T(0)K_{\pm}(n,n)}.
\ee

H.\ref{hypo scat q} (iv) and (\ref{W q}) 
imply the following behavior of $\hat h_\mp(\lambda, n)$ as 
$\lambda \rightarrow \rho_j$
\bea       \label{h(rho) 2}      \nonumber
\lim_{\lambda \rightarrow \rho_j} \hat h_\mp(\lambda, n) 
&=& \pm\gamma_{\pm,j} \hat \psi_\pm(\rho_j,n) 
\lim_{\lambda \rightarrow \rho_j} 
\frac{W_{n-1}(\hat \psi_{q,\pm}(\rho_j), \hat \psi_{q,\mp}(\lambda))}
{(\lambda - \rho_j)T(\lambda)} \\  \label{h(rho)}
&=& \gamma_{\pm,j} \hat \psi_\pm(\rho_j, n)
\big(\res_{\rho_j} T(\lambda)\big)^{-1}\prod_{l=0}^{2g+1}\sqrt{\rho_j - E_l},
\eea
where $\hat h_\pm$ are defined as in (\ref{def eigen}).

By virtue of the consistency condition $T(w)\overline{R_+(w)} = -
\overline{T(w)}R_-(w)$ we obtain
\bea              \nonumber
\lefteqn{\overline{h_{\pm}(w, n)} + R_{\pm}(w)h_{\pm}(w, n) =}    \\ \nonumber
&=& \frac{1}{\overline{T(w)}} 
\Big( \psi_{\mp}(w, n) + \overline{R_{\mp}(w)} 
\overline{\psi_{\mp}(w, n)} \Big) + \frac{R_{\pm}(w)}{T(w)} 
\Big( \overline{\psi_{\mp}(w, n)} 
+ R_{\mp}(w) \psi_{\mp}(w, n) \Big) \\ \nonumber
&=& \psi_{\mp}(w, n) \Big( \frac{1}{\overline{T(w)}} 
+ \frac{R_{\pm}(w)R_{\mp}(w)}{T(w)}\Big)
+ \overline{\psi_{\mp}(w, n)}\Big(\frac{\overline{R_{\mp}(w)}}
{\overline{T(w)}} + \frac{R_{\pm}(w)}{T(w)} \Big)\\   \nonumber          
&=& \psi_{\mp}(w, n) T(w), \qquad \qquad |w|=1.
\eea 
If we eliminate $R_{\mp}(w)$ from the last equation
and (\ref{h3 q}) we see
\bea                   \nonumber
T(w)R_{2g+2}^{-1/2}(w)\big(\hat \psi_+(w,n) \hat \psi_-(w,n) 
- \hat h_+(w,n) \hat h_-(w,n)\big) 
&\!=\!& \\       \label{def G}  
\frac{\prod_j (\lam(w)-\mu_j)}{R_{2g+2}^{1/2}(w)}
\big(\overline{h_\pm(w,n)} \psi_\pm(w,n) 
- \overline{\psi_\pm(w,n)} h_\pm(w,n)\big)   
&\!=:\!& G(w,n),
\eea
for $|w|=1$.
Observe that $G(\overline{w}, n) = \overline {G(w, n)} = G(w, n)$, 
$|w|=1$, since $\overline{\hat h_\pm} \hat \psi_\pm 
- \overline{\hat \psi_\pm} \hat h_\pm$ and $R_{2g+2}^{-1/2}(w)$ are odd 
functions for $|w|=1$.
The function $G(w, n)$ can be continued analytically on $\W$ since 
the difference $\hat \psi_+ \hat \psi_- - \hat h_+ \hat h_-$ vanishes at the 
poles $w(\rho_j)$ of $T(w)$ by (\ref{h(rho)}). 
Note that the product $\hat \psi_+ \hat \psi_-$ and hence also $\hat h_+ \hat h_-$
do not have poles at $w(\mu_j)$. Moreover, since $\W$ is just the image of
the upper sheet, we can extend it to a compact Riemann surface $\ti{\W}$ by
adding the image of the lower sheet. Now by $G(\overline{w}, n) = G(w, n)$
we can extend $G$ to $\ti{\W}$ by setting $G(w, n)= G(w^{-1},n)$ for $|w|>1$.

Now let us investigate the behavior at the band edges:
If $w_l\neq w(\mu_j)$, we obtain by (\ref{lim T}), (\ref{h3 q}),
and real-valuedness of $\hat{\psi}_\pm$ at the band edges that 
\bea           \nonumber
\lefteqn{\lim_{w \rightarrow w_l}  R_{2g+2}^{1/2}(w)  \prod_j(\lam(w)-\mu_j)
h_\mp(w, n) \overline{\psi_\mp(w,n)}}\\ \nonumber
&=& 
\lim_{w \rightarrow w_l} \frac{R_{2g+2}^{1/2} \prod_j(\lam-\mu_j)}{T}
\Big( \overline{\psi_\pm} + R_\pm \psi_\pm \Big) \overline{\psi_\mp}
\\                              \nonumber
&=& \lim_{w \rightarrow w_l} \frac{R_{2g+2}^{1/2} \prod_j(\lam-\mu_j)}{T}
\Big( (R_\pm + 1) \psi_\pm + \overline{\psi_\pm} - \psi_\pm \Big) 
\overline{\psi_\mp} = 0. 
\eea 
If $w_l = w(\mu_j)$, the same calculation shows that
\begin{align*}
&\lim_{w \rightarrow w_l}  R_{2g+2}^{1/2}(w)  \prod_j(\lam(w)-\mu_j)
h_\pm(w,n) \overline{\psi_\pm(w,n)} \\
&= (-1)^{l+1} C_+(n)C_-(n) \lim_{w \rightarrow w_l}
R_{2g+2}^{1/2}(w)\frac{R_\pm(w) - 1}{T(w)} = 0 
\end{align*}
by (\ref{lim T}), where we used 
$\psi_\pm(w,n)= \I^l C_\pm(n)(\lambda(w)-\mu_j)^{-1/2} + O(1)$.

Consequently
$R_{2g+2}(w) G(w,n)$ is continuous at $w=w_l$ and vanishes at the band edges.
Thus the singularities of $R_{2g+2}^{1/2}(w) G(w,n)$ at $w_l$ are removable.
Furthermore, $R_{2g+2}^{1/2}(w) G(w,n)$ is purely imaginary for $|w|=1$ and
real on the slits and hence must vanish at $w_l$ by continuity.
So the singularities of $G(w,n)$ at $w_l$ are removable as well.
Thus $G$ is holomorphic on all of $\ti{\W}$ and vanishes at $w=0$, that is,
$G(w,n) \equiv 0$ which implies (compare (B4))
\bea   \nonumber     
\lefteqn{\lim_{w \rightarrow 0}
\Big(\psi_+(w,n) \psi_-(w,n)- h_+(w,n)h_-(w,n)\Big)} \\ \nonumber 
&=& K_+(n,n)K_-(n,n) - (T(0)^2 K_+(n,n)K_-(n,n))^{-1} = 0. 
\eea 
Using (\ref{a+a- q}) we finally obtain
from  
$T(0)^2 = \big(K_+(n, n) K_-(n, n)\big)^{-2}$ that
\be
a_+(n) = a_-(n) \equiv a(n), \qquad \forall n \in \Z.
\ee
It remains to prove $b_+(n)=b_-(n)$. Proceeding as for $G(w,n)$ we can show
that
\bea                   \nonumber
&&T(w)R_{2g+2}^{-1/2}(w)\big(\hat \psi_+(w,n) \hat \psi_-(w,n+1) 
- \hat h_+(w,n+1) \hat h_-(w,n)\big)=\\
&&\quad \frac{\prod_j (\lam(w)-\mu_j)}{R_{2g+2}^{1/2}(w)}
\big(\overline{h_+(w,n+1)} \psi_+(w,n) 
- \overline{\psi_+(w,n)} h_+(w,n+1)\big)
\eea
is a constant equal to $-1/a(n)$. Thus
\bea \nn
W(w,n) &:=& a(n) \left(\psi_+(w,n) \psi_-(w,n+1) - h_+(w,n+1) h_-(w,n) \right)\\
&=& -\frac{R_{2g+2}^{1/2}(w)}{T(w) \prod_j(\lam(w)-\mu_j)}.
\eea
and computing the asymptotics at $w=0$ (compare (\ref{B4jost})) we see
\be
0= W(w,n) - W(w,n-1) = A_+(0)A_-(0) (b_+(n) - b_-(n))
\ee
and in particular $b_+(n) = b_-(n) \equiv b(n)$.

Our operator $H$ has the correct norming constants since as in (\ref{derivative alpha}) it follows
\be
\sum_{n \in \Z} \hat \psi_+(\rho_j, n) \hat \psi_-(\rho_j, n) =
\big(\res_{\rho_j} T(\lambda)\big)^{-1} \prod_{l=0}^{2g+1}\sqrt{\rho_j - E_l} 
\ee 
and by (\ref{h(rho) 2}),
\[
\sum_{n \in \Z} \hat \psi_{\pm}(\rho_j, n) \hat \psi_{\pm}(\rho_j, n) = \gamma_{\pm,j}^{-1}.
\]
\end{proof}

\section*{Acknowledgments}
I.E.\ thanks A.\ Boutet de Monvel for the kind hospitality of
University Paris-7, where part of this work was done.
G.T.\ thanks Peter Yuditskii for several helpful discussions and hints with respect to
the literature. 
We thank Mark Losik for help with respect to literature.


\begin{thebibliography}{XXXX}
\bibitem{baeg}
J. Bazargan and I. Egorova, {\em Jacobi operator with step-like asymptotically periodic
coefficients}, Mat. Fiz. Anal. Geom. {\bf 10}, no. 3, 425--442 (2003).
\bibitem{bdme} A. Boutet de Monvel and I. Egorova, {\em Transformation operator for
Jacobi matrices with asymptotically periodic coefficients}, J. Difference Eqs. Appl. {\bf 10}, 711-727 (2004).
\bibitem{bdme2} A. Boutet de Monvel and I. Egorova, {\em  The Toda lattice with step-like initial data. Soliton asymptotics}, Inverse Problems {\bf 16}, no. 4, 955--977 (2000).
\bibitem{bght} W. Bulla, F. Gesztesy, H. Holden, and G. Teschl {\em
Algebro-Geometric Quasi-Periodic Finite-Gap Solutions of the Toda and Kac-van
Moerbeke Hierarchies}, Memoirs of the Amer. Math. Soc. {\bf 135/641}, (1998).
\bibitem{dinv4} K. M. Case, {\em Orthogonal polynomials from the viewpoint of
scattering theory}, J. Math. Phys. {\bf 14}, 2166--2175 (1973).
\bibitem{dinv5} K. M. Case, {\em The discrete inverse scattering problem in one
dimension}, J. Math. Phys. {\bf 15}, 143--146 (1974).
\bibitem{caopt} K. M. Case, {\em Orthogonal polynomials II},  J. Math. Phys.
{\bf 16}, 1435--1440 (1975).
\bibitem{dinv2} K. M. Case, {\em On discrete inverse
scattering problems. II}, J. Math. Phys. {\bf 14}, 916--920 (1973).
\bibitem{dinv3} K. M. Case and S. C. Chiu {\em The discrete version of the
Marchenko equations in the inverse scattering problem}, J. Math. Phys. {\bf 14},
1643--1647 (1973).
\bibitem{dinv1} K. M. Case and M. Kac, {\em A discrete version of the inverse
scattering problem}, J. Math. Phys. {\bf 14}, 594--603 (1973).
\bibitem{fad} L. Faddeev and L. Takhtajan, {\em Hamiltonian Methods in the
Theory of Solitons}, Springer, Berlin, 1987.
\bibitem{coj} P. A. Cojuhari, {\em Finiteness of the discrete spectrum of Jacobi
matrices} (Russian), Investigations in differential equations and mathematical analysis
{\bf 173}, "Shtiintsa", Kishinev, 80--93 (1988).
\bibitem{fir} N.E. Firsova, 
{\em The direct and inverse scattering problems for the one-dimensional perturbed Hill operator},
Math. USSR, Sb. {\bf 58}, 351--388 (1987).
\bibitem{fl2} H. Flaschka, {\em On the Toda lattice. II}, Progr. Theoret. Phys.
{\bf 51}, 703--716 (1974).\bibitem{ggkm} C. S. Gardner, J. M. Green, M. D. Kruskal, and R. M. Miura,
{\em A method for solving the Korteweg-de Vries equation}, Phys. Rev.
Letters {\bf 19}, 1095--1097 (1967).
\bibitem{gnp} F. Gesztesy, R. Nowell and W. P\"otz, {\em One-dimensional scattering theory for quantum
systems with nontrivial spatial asymptotics},
Differ. Integral Equ. {\bf 10}, No.3, 521--546 (1997).
\bibitem{gerass} J. S. Geronimo and W. Van Assche, {\em Orthogonal polynomials
with asymptotically periodic recurrence coefficients},
J. App. Th. {\bf 46}, 251--283 (1986).
\bibitem{gu} G.S. Guseinov, {\em The inverse problem of scattering theory for a
second-order difference equation on the whole axis}, Soviet Math. Dokl., {\bf 17},
1684--1688 (1976).
\bibitem{gu2} G.S. Guseinov, {\em The determination of an infinite Jacobi matrix
from the scattering data}, Soviet Math. Dokl., {\bf 17}, 596--600 (1976).
\bibitem{gu3} G.S. Guseinov, {\em Scattering problem for the
infinite Jacobi matrix}, Izv. Akad. Nauk Arm. SSR, Mat. {\bf 12},
365--379 (1977).
\bibitem{kumi} E.A. Kuznetsov and A.V. Mikha\u\i lov, {\em Stability of stationary waves in
nonlinear weakly dispersive media}, 
 \v Z. \`Eksper. Teoret. Fiz. {\bf 67} no. 5, 1717--1727 (1974). (in
Russian)
\bibitem{mar} V.A. Marchenko, {\em Sturm--Liouville Operators and Applications}, 
Birkh\"auser, Basel, 1986.
\bibitem{par} T.\ Parthasarathy, {\em On Global Univalence Theorems}, LNM
577, Sprin\-ger, Berlin, 1983.
\bibitem{perco} L.\ Percolab, {\em The inverse problem for the periodic
Jacobi matrix}, Theor.\ funk., funk.\ an., pril.\ {\bf 42}, 107--121
(1984), in Russian.
\bibitem{tosc} G. Teschl, {\em Oscillation theory and renormalized oscillation
theory for Jacobi operators}, J. Diff. Eqs. {\bf 129}, 532--558 (1996).
\bibitem{tist} G. Teschl, {\em Inverse scattering transform for the Toda hierarchy},
Math. Nach. {\bf 202}, 163--171 (1999).
\bibitem{tivp} G. Teschl, {\em On the initial value problem for the Toda and
Kac-van Moerbeke hierarchies}, AMS/IP Studies in Advanced Mathematics {\bf 16}, 375--384 (2000).
\bibitem{tjac} G. Teschl, {\em Jacobi Operators and Completely Integrable Nonlinear Lattices}, Math. Surv. and Mon. {\bf 72}, Amer. Math. Soc., Rhode Island, 2000.
\bibitem{ta} M. Toda, {\em Theory of Nonlinear Lattices}, 2nd enl. ed.,
Springer, Berlin, 1989.
\bibitem{tsu} M. Tsuji, {\em Potential Theory in modern Functional Analysis}, Maruzen, Tokyo, 1959.
\bibitem{voyu} A. Volberg and P. Yuditskii, {\em On the inverse scattering problem for Jacobi Matrices with the Spectrum on an Interval, a finite systems of intervals or a Cantor set of positive length},
Commun. Math. Phys. {\bf 226},  567--605 (2002).
\bibitem{voza} V. Voichick and L. Zalcman, {\em Inner and outer functions on Riemann surfaces},
Proc. Amer. Math. Soc. {\bf 16}, 1200-1204 (1965).
\end{thebibliography}
\end{document}